\documentclass[10pt]{amsart}

\usepackage{epic,eepic}
                                                                               
\usepackage{latexsym}
\usepackage{amssymb}
\usepackage{graphicx}
                                                                                
\usepackage{latexsym}
\usepackage{amssymb}
                                                                                
\numberwithin{equation}{section}
\newtheorem{thm}{Theorem}[section]
\newtheorem{lemma}[thm]{Lemma}

\newtheorem{prop}[thm]{Proposition}
\newtheorem{conj}[thm]{Conjecture}

\newtheorem{theorem}{Theorem}
\newenvironment{remark}{\begin{trivlist}\item {\bf
        Remark.\,}}{\mbox{}\hfill$\square$\hfill\end{trivlist}}

\newenvironment{remarks}{\begin{trivlist}\item {\bf
        Remarks.\,}}{\mbox{}\hfill$\square$\end{trivlist}}
\newcommand{\excise}[1]{}

    {\begin{eqnarray}}
    {\end{eqnarray}$\!\!$}

\newenvironment{eq*}%
    {\begin{eqnarray*}}
    {\end{eqnarray*}$\!\!$}
    {\begin{equation}}
    {\end{equation}$\!\!$}
                                                                                
\newenvironment{eqn*}%
    {\begin{equation*}}
    {\end{equation*}$\!\!$}
\def\gln{\mathsf{GL_n}}
\def\sln{\mathsf{SL_n}}
\def\pp{{\bf p}}

\def\CC{{\mathbb C}}
\def\Gr{{\bf Gr}}
\def\MV{{M\!}}

\def\cplx{\CC} 
\newcommand{\cdim}{\mbox{dim}} 
\newcommand{\cspan}{\mbox{span}} 
 
\newcommand{\kost}{\mathbb{K}} 
\newcommand{\base}{\eta} 
\newcommand{\reals}{\mathbb{R}} 
\newcommand{\polytinv}{\mathbb{C}[t,t^{-1}]}

\def\vs{{\mathfrak t}_\reals} 
\def\momm{\Phi} 
\def\qq{{\bf q}}
\def\torus{T}
\def\det{\mbox{det}}
\def\poly{\mathcal P}
\def\nel{\mathbf n}

\def\fone{\chi_1}
\def\ftwo{\chi_2}
\def\bz{{\mathbf 0}}                                                                                
\def\diredge{\rightarrow}
\def\gp{\mathsf G}
\def\torus{\mathsf T}
\def\unip{\mathsf N}
\def\borel{\mathsf B}
\def\calg{{\mathcal C}_\gp}
\def\palg{{\CC[\unip]}}
\def\th{^{\rm th}\!}

\def\convolution{{convolution }}
\def\BD{{Beilinson-Drinfeld }}

\def\conv{{*}}
\def\convC{{*_{\mathbb C}}}
\def\qq{{\bf q}}
\def\rr{{\bf r}}
\def\ss{{\bf s}}
\def\xx{{\bf x}}
\def\KK{{\mathbb K}}
\def\ZZ{{\mathbb Z}}
\def\PP{{\mathbb P}}

\def\isom{{\mathcal I}}
\def\cluster{{\bf c}}
\def\OO{{\mathcal O}}
\def\KKK{{\mathcal K}}
\def\RR{{\mathbb R}}

\newcommand{\yz}{Y_1} 
\newcommand{\ys}{Y_2}

\begin{document}

\title[The algebra of Mirkovi\'c-Vilonen cycles in type A]
    {The algebra of Mirkovi\'c-Vilonen cycles in type A} \author{Jared 
Anderson} 
\address{Institute for Advanced Study\\Princeton, NJ\\USA} 
\email{jareda@ias.edu}
                                                                                
\author{Mikhail Kogan} 
\thanks{Both authors are  supported by the National Science Foundation under agreement No. 
DMS-0111298} 
\address{Institute for Advanced Study\\Princeton, NJ\\USA} 
\email{mish@ias.edu}

\begin{abstract} 
\noindent 
Let $\Gr$ be the affine Grassmannian for 
a connected complex reductive group $\gp$\@.  
Let~$\calg$ be the complex vector space spanned by  
(equivalence classes of) Mirkovi\'c-Vilonen cycles in $\Gr$\@.  
The Beilinson-Drinfeld Grassmannian can be used 
to define a convolution product on MV-cycles, making 
$\calg$ into a commutative algebra.  
We show, in type~A, that $\calg$ is isomorphic to 
$\palg$, the algebra of functions on the unipotent radical $\unip$ 
of a Borel subgroup of $\gp$; 
then each MV-cycle defines a polynomial in $\palg$, 
which we call an MV-polynomial.  
We conjecture that those MV-polynomials which 
are cluster monomials for a 
Fomin-Zelevinsky 
cluster algebra structure on $\palg$ 
are naturally expressible as determinants, 
and we conjecture a formula for many of them.

\vspace{3mm}
                                                                                                                 
\noindent (Mathematics subject classification number: 14L35)
                                                                                                                 
\end{abstract}
                                                                                                                 
\maketitle

\vspace{-3mm}
                                                                                                                 
\begin{center}{ \em  This paper is dedicated to Robert MacPherson
on the occasion of his 60th birthday.}
\end{center}

\vspace{3mm}
                                                                                                                 
\section{Introduction} 
\label{sec:intro}

Since Mirkovi\'c and Vilonen first announced their discovery 
of a new geometric canonical basis in representation theory~\cite{MV1}, 
the study of their algebraic varieties---here called MV-cycles---has 
grown quickly, with applications to geometry, 
representation theory, and combinatorics.  
Let us begin by mentioning a few of these.

In the extended version of their 
original paper, Mirkovi\'c and Vilonen~\cite{MV2} used MV-cycles to 
study the relation between the representation theory of a 
complex reductive group and the equivariant perverse 
sheaves on the affine Grassmannian $\Gr$ for the Langlands dual 
group.  
The MV-cycles are subvarieties of $\Gr$ giving a canonical basis 
for every irreducible representation.
Braverman and Gaitsgory~\cite{BG} used MV-cycles to give a geometric definition of
crystal graphs.  Gaussent and Littelmann~\cite{GL} gave an interpretation of
Littelmann's path model using MV-cycles.  Kamnitzer~\cite{Kamnitzer}
discovered the relation between the combinatorics of the MV-polytopes (moment
map images of MV-cycles) and the Berenstein-Zelevinsky combinatorics 
used to compute Littlewood-Richardson coefficients~\cite{BZ}.  

Essential to the Mirkovi\'c-Vilonen results 
was the Beilinson-Drinfeld Grassmannian. This is a space which allows one 
to define a convolution product of MV-cycles (sometimes called a fusion product), 
making the vector space spanned by all MV-cycles into a commutative algebra.  
The purpose of this paper is to discuss this algebra and its combinatorial 
properties in type~A.  In particular, we identify it with 
a polynomial algebra, and conjecture a combinatorial formula for 
some of the polynomials.

We will continue in the spirit of~\cite{AK}, in which the 
authors computed all MV-polytopes in type~A; we use mostly elementary tools 
and combinatorial models, and we avoid abstraction, always preferring to make 
concrete choices that make it easy to work with examples.  
There are two reasons for working in this way: 
(1) Our results and conjectures can be stated in this elementary language.  
(2) The literature about MV-cycles often gives the inaccurate impression 
that considerable heavy machinery and abstruse mathematics is required to 
understand them.

Let us now discuss the main results and conjectures of the paper. As
mentioned above, the MV-cycles are certain algebraic varieties lying inside
the affine Grassmannian $\Gr$ for a connected complex reductive group $\gp$\@. 
The standard short definition of $\Gr$ is that $\Gr=\gp(\KKK)/\gp(\OO)$, 
where $\KKK=\CC((t))$ and 
$\OO=\CC[[t]]$\@.  But, for our purposes, this definition is difficult 
to work with and provides little geometric intuition; we prefer to use a version 
of Lusztig's lattice model~\cite{L} for $\Gr$ in type~A, which we will 
discuss in 
Section~\ref{sec:lattice-model}.  Although the affine Grassmannian is an 
infinite-dimensional space, it 
is an increasing union of finite-dimensional pieces (formally an ind-scheme).  
In any example we will be interested in, we may work completely 
within one of these pieces.  

One can also give a concise definition of MV-cycles: 
A maximal torus $\torus$ of $\gp$ acts on $\Gr$ and has fixed point set
canonically isomorphic to the coweight lattice $\Lambda={\rm Hom}
(\CC^*,T)$:  $\lambda\in\Lambda \leftrightarrow \underline{\lambda}\in\Gr$\@. 
Moreover, this action is Hamiltonian and has a naturally defined
moment map $\momm:\Gr\to \vs={\rm Lie}(\torus)$\@. Consider the Morse flow of a generic
component $\xi :\Gr\to \RR$ of the moment map. Let $S_\lambda $ be the
space of points of $\Gr$ which converge to $\underline{\lambda}$\@. 
Similarly, define $T_\lambda$ for the
inverse flow. (Actually, later we will pick a concrete component
$\xi$\@.)  Then MV-cycles are defined to be the irreducible components of the
closures of $S_\lambda\cap T_\mu$\@.  Their moment map images are called MV-polytopes.
In Section~\ref{sec:Definition of MV-cycles} 
we will repeat this definition in the language of the lattice model, 
and also recall from~\cite{AK} how to describe and parametrize the 
irreducible components~(Theorem~\ref{thm:param}), a description we will use heavily throughout.

Our primary concern is with the convolution product of MV-cycles:  
given two MV-cycles~$M$ and~$N$, this expresses their product as a linear
combination of MV-cycles
\begin{equation}
\label{introconv}
M\conv N=\sum c_{MN}^K K.
\end{equation} 
This convolution is defined using the Beilinson-Drinfeld Grassmannian.  
This space is a fibration over $\cplx$ whose fiber over zero is
$\Gr$ and whose fibers over $\cplx^*$ are isomorphic to $\Gr\times \Gr$\@. 
Within this family, it is possible to construct a family $M\convC N$ which
degenerates $M\times N$ to the union of $K$'s with multiplicities
$c_{MN}^K$\@. We will present in Section~\ref{sec:conv} 
a lattice model for this family, which will
allow us to identify, in type~A, 
the leading term on the right side of~(\ref{introconv}).   
(Propositions~\ref{prop:order} and~\ref{prop:lead}).

This will result in an algebra $\calg$ for which the MV-cycles provide 
a basis, as a complex vector space.  
(Actually, we will have some simple equivalence relation
on MV-cycles, and $\calg$ will be generated by the equivalence classes.)  
It was conjectured by the first author~\cite{Anderson-thesis}
that $\calg$ is isomorphic as a Hopf algebra to the ring of functions $\palg$ 
on the unipotent radical $\unip$ of a Borel subgroup of $\gp$\@.  
In Section~\ref{sec:convolution algebra} we will prove part of this conjecture 
in type~A, where $\unip$ is the group of upper triangular unipotent 
matrices.  Using the properties we have proved about convolution, 
we will construct a map $\isom: \palg \rightarrow \calg$ and 
prove that it is an algebra isomorphism (Theorem~\ref{thm}).  
So corresponding to each MV-cycle $M \in \calg$ 
we have a polynomial $\poly \in \palg$; 
we call these MV-polynomials.

Section~\ref{sec:determinant} addresses the computation of MV-polynomials and their 
relation to cluster algebras. 
The simplest MV-cycles are Richardson varieties (intersections of opposite 
Schubert varieties) and the corresponding MV-polynomials are exactly the 
minors in $\palg$\@. But understanding MV-polynomials in general 
seems to be a difficult problem and we will only conjecture how to compute some of them. 

Of crucial theoretical and practical importance to our conjectures
is the Fomin-Zelevinsky theory of cluster algebras~\cite{FZ1,FZ2,BFZ}.
Part of this structure is a canonical set of polynomials in $\palg$---the {\em cluster 
monomials}. We conjecture (Conjecture~\ref{conj:cluster}) that these are themselves
MV-polynomials, and, moreover,
are all naturally expressible as certain determinants.

Although the inductive computation of MV-polynomials
from the definition of convolution of MV-cycles is difficult,
the computation of the cluster
monomials is straightforward and has been programmed onto a computer.  
(The algorithm is specified in Appendix~\ref{appendix} and our data for 
type~$\rm A_5$ may 
be found at~\cite{AK*}.)  
This provided many examples to work with, 
based upon which we sharpened Conjecture~\ref{conj:cluster} to Conjecture~\ref{conj:zeros}, 
which gives an explicit determinantal formula for a large class of MV-polynomials.  

Finally, Conjecture~\ref{conj:samefunction} relates the MV-polynomials that are cluster 
monomials to the geometry of MV-cycles. 
The most intriguing thing about all our conjectures is the connection between cluster 
algebras and MV-cycles, which we hope will eventually be understood on a deeper level.

To summarize, 
Section~\ref{sec:MV} recalls what we need from~\cite{AK}: 
the definition of the lattice model for~$\Gr$ and the 
explicit parametrization and description of MV-cycles in type~A.  
Section~\ref{sec:conv} discusses the convolution product of 
MV-cycles using a lattice model for the Beilinson-Drinfeld Grassmannian.  
In Section~\ref{sec:convolution algebra} we prove 
that the convolution algebra of MV-cycles is isomorphic to 
the polynomial algebra $\palg$, and use this isomorphism 
to define MV-polynomials.  
Section~\ref{sec:determinant} contains conjectural determinantal 
formulas for some MV-polynomials.   Appendix~A is a terse 
summary of what we use from cluster algebras.

\smallskip
                                                                                
\noindent{\bf Acknowledgments.} Both authors thank J.~Kamnitzer, A.~Knutson, 
B.~Leclerc, R.~MacPherson,
I.~Mirkovi\'c, D.~Nadler and A.~Zelevinsky for helpful discussions and suggestions.

\section{MV-Cycles}
\label{sec:MV}

In this section, we review the type A description of MV-cycles 
and their parametrization by the Kostant parameter set, as set forth in~\cite{AK}.  
We use throughout the lattice model of the affine Grassmannian, 
the space in which the MV-cycles live.

\subsection{Lattice model} 
\label{sec:lattice-model}
                                                                                
Following Lusztig~\cite{L}, we will introduce the lattice model of the affine 
Grassmannian~$\Gr$ for 
$\gp=\gln$\@.  Points in $\Gr$ are subspaces of a certain infinite-dimensional complex vector 
space~$X$, 
satisfying some extra conditions.  The vector space $X$ is defined by specifying a basis: if 
$e_1, e_2, \ldots, e_n$ denotes the standard basis for~$\cplx^n$, then $X$ is the span of 
$t^je_i$ for $1\leq i\leq n$, $-\infty<j<\infty$, where we regard this as a symbol with two 
indices.  We usually picture these basis vectors in an array consisting of $n$ columns, 
infinite in both directions. Let $t$ be the invertible linear operator on $X$ that sends 
each~$t^je_i$ to $t^{j+1}e_i$\@.  Then $\Gr$ consists of those subspaces $Y$ of $X$ such that 
\begin{enumerate} \item $Y$ is closed under the action of $t$\@. \item $t^NX_0\subseteq 
Y\subseteq t^{-N}X_0$ for some $N$,
  where $X_0$ is the span of those $t^je_i$ with $j\geq0$\@. \end{enumerate} We call such $Y$ 
{\em lattices}.

\begin{figure}[ht] \begin{picture}(400,73)
                                                                                
\put(0,10){ \setlength{\unitlength}{.5mm}
                                                                                
\def\smalldot{{$\bullet$}} \def\tinydot{{\tiny$\bullet$}}
                                                                                
\put(0,0){$t^2$} \put(0,10){$t$} \put(0,20){$1$} \put(0,30){$t^{-1}$} \put(0,40){$t^{-2}$} 
\put(0,-10){$\vdots$}
                                                                                
\put(20,0){\smalldot} \put(30,0){\smalldot} \put(40,0){\smalldot} \put(50,0){\smalldot} 
\put(60,0){\smalldot} \put(70,0){\smalldot}
                                                                                
\put(20,10){\smalldot} \put(30,10){\smalldot} \put(40,10){\smalldot}
                                                                                
\put(20,20){\smalldot} \put(20,30){\smalldot} \put(20,40){\smalldot}
                                                                                
\put(45,-5){\tinydot} \put(45,-9){\tinydot} \put(45,-13){\tinydot}
                                                                                
\put(100,0){$t^2$} \put(100,10){$t$} \put(100,20){$1$} \put(100,30){$t^{-1}$} 
\put(100,40){$t^{-2}$} \put(100,-10){$\vdots$}
                                                                                
\put(120,0){\smalldot} \put(130,0){\smalldot} \put(140,0){\smalldot} \put(150,0){\smalldot} 
\put(160,0){\smalldot} \put(170,0){\smalldot}
                                                                                
\put(120,10){\smalldot} \put(130,10){\smalldot} \put(140,10){\smalldot} 
\put(160,10){\smalldot} \put(170,10){\smalldot}
                                                                                
\put(120,20){\smalldot} \put(130,20){\smalldot} \put(140,20){\smalldot}
                                                                                
\put(120,30){\smalldot} \put(120,40){\smalldot}
                                                                                
                                                                                
\put(132,22.5){\line(1,0){10}} \put(132,12.5){\line(1,0){10}} \put(132,22){\line(-1,2){10}} 
\put(143,14.5){\line(1,0){20}} \put(152,2){\line(1,1){10}} \put(172,12){\line(-1,0){10}}
                                                                                
                                                                                
\put(141.5,12){\smalldot} \put(161.5,12){\smalldot}
                                                                                
                                                                                
{\tiny
                                                                                
\put(118,39){$1$} \put(128,18){$3$} \put(138,18){$4$}
                                                                                
\put(128,8){$3$} \put(138,8){$4$}
                                                                                
\put(152,-2){$1$} \put(162,8){$2$} \put(172,8){$-1$}
                                                                                
\put(142,15){$-1$} \put(162,15){$6$}}
                                                                                
\put(145,-5){\tinydot} \put(145,-9){\tinydot} \put(145,-13){\tinydot}
                                                                                
                                                                                
\put(200,0){$t^2$} \put(200,10){$t$} \put(200,20){$1$} \put(200,30){$t^{-1}$} 
\put(200,40){$t^{-2}$} \put(200,-10){$\vdots$}
                                                                                
\put(220,0){\smalldot} \put(230,0){\smalldot} \put(240,0){\smalldot} \put(250,0){\smalldot} 
\put(260,0){\smalldot} \put(270,0){\smalldot}
                                                                                
 \put(230,10){\smalldot} \put(240,10){\smalldot} \put(250,10){\smalldot} 
\put(260,10){\smalldot} \put(270,10){\smalldot}
                                                                                
\put(230,20){\smalldot} \put(240,20){\smalldot} \put(250,20){\smalldot} 
\put(260,20){\smalldot} \put(270,20){\smalldot} \put(270,30){\smalldot} 
\put(270,40){\smalldot}
                                                                                
                                                                                
\put(232,22.5){\line(1,0){30}} \put(232,12.5){\line(1,0){20}} \put(232,14.5){\line(1,0){28}} 
\put(232,2.5){\line(1,0){20}}
                                                                                
\put(222,2){\line(1,1){10}}

                                                                                
\put(262,12){\line(1,2){10}} \put(262,22){\line(1,2){10}}
                                                                                
\put(231.2,2){\line(0,1){10}}
                                                                                

\put(231.5,12.5){\smalldot} \put(241.5,12.5){\smalldot} 
\put(251.5,12.5){\smalldot} \put(258.5,12.5){\smalldot}
                                                                                
\put(230,2){\smalldot}
                                                                                
                                                                                
{\tiny
                                                                                
\put(232,-2){$1$}\put(242,-2){$2$}\put(252,-2){$3$}
                                                                                
\put(232,8.5){$1$} \put(242,8){$2$} \put(252,8){$3$} \put(232,4){$1$} \put(222,-2){$1$}
                                                                                
\put(230.5,15.5){$1$}\put(240.5,15.5){$2$}\put(250.5,15.5){$3$} \put(257.5,15.5){$1$}
                                                                                

\put(229,23){$1$} \put(239,23){$2$} \put(249,23){$3$} \put(259,23){$2$} 
\put(269,43){$1$}
                                                                                
\put(262,8){$1$} \put(272,28){$1$} \put(245,-5){\tinydot} \put(245,-9){\tinydot} 
\put(245,-13){\tinydot} }

}
                                                                                
\end{picture} 
\caption{Examples of lattices.} 
\label{fig:lattice-examples} 
\end{figure}
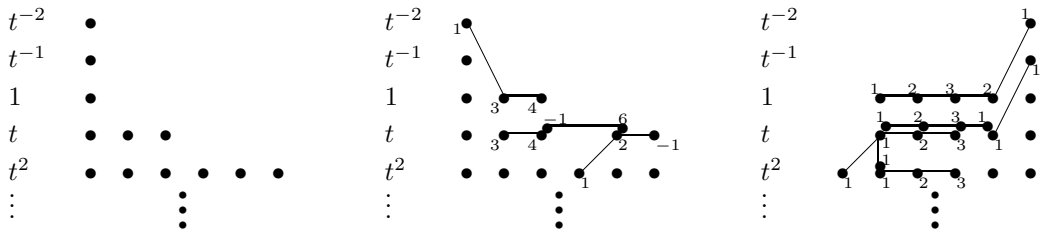

Three examples of lattices are illustrated by the pictures in 
Figure~\ref{fig:lattice-examples}. Each dot, or set of dots connected by line segments, 
represents a vector; these, together with all the dots $t^je_i$, $j\geq 3$ below the pictures 
are a $\CC$-basis for the lattice. The first picture represents the span of the basis vectors 
drawn, namely $t^je_1$ with $j\geq -2$, $t^je_2$, $t^je_3$ with $j\geq 1$, and $t^je_4$, 
$t^je_5$, $t^je_6$ with 
$j\geq 2$\@.  The second picture represents the lattice spanned by 
the vectors $t^{-2}e_1+3e_2+4e_3$, $3te_2+4te_3$, $-te_3+6te_5$, $t^2e_4+2te_5-te_6$, 
and all the vectors of the form $t^je_i$ corresponding to the unconnected dots.  
Similarly for the third picture.

Lattices spanned by vectors $t^je_i$ are in one-to-one correspondence with the 
coweights of 
$\gln$ (which we will regard as n-tuples of integers).  The first picture in 
Figure~\ref{fig:lattice-examples} corresponds to the coweight $(2,-1,-1,-2,-2,-2)$ according 
to our convention.  A coweight with nonincreasing entries, like this one, is said to be {\em 
dominant}; one with nondecreasing entries is said to be {\em antidominant}.  
In general, if 
$\lambda$ is a coweight then we denote the corresponding lattice by $\underline{\lambda}$\@.  
Also, any coweight $\lambda=(\lambda_1,\lambda_2,\dots,\lambda_n)$ acts on $X$ by 
shifting the $i\th$ column by $\lambda_i$: $\lambda \cdot t^je_i = t^{j-\lambda_i}e_i$\@. 
This allows us to form the lattice $\lambda \cdot Y$, the {\em shift} of $Y$ by $\lambda$, 
by setting $\lambda \cdot Y = \{ \lambda \cdot v \mid v\in Y \}$\@. It is then clear that 
$\lambda\cdot\underline \mu=\underline{\lambda+\mu}$\@.


Note that if $N$ is fixed in condition (2), then the set of all lattices
is a finite-dimensional projective variety  $\Gr^N\!$\@.
So $\Gr$ is an increasing union of complex projective varieties.  
Moreover,
although the definition of $\Gr$ at first sounds very
infinite-dimensional, in practice everything is finite-dimensional: all
but a finite-dimensional subspace of each lattice is spanned by basis
vectors, and the varieties we will study will each be contained in one of
the $\Gr^N\!$\@.

The affine Grassmannian $\Gr$ for $\gln$ is not connected. Each connected component contains 
all lattices~$Y$ with fixed dimension $\dim_0 Y= \dim(Y/Y\cap X_0)-\dim (X_0/Y\cap X_0)$ 
relative to $X_0$\@.
All the connected components are isomorphic to each other by 
shifts.  Also, the connected component where $\dim_0 Y=0$ is the affine 
Grassmannian for $\sln$\@.

\subsection{Torus action and moment map} 
\label{sec:torus action and moment map}

There is an algebraic action of an $(n-1)$-dimensional torus~$\torus$ on $\Gr$:  $(\CC^*)^n$ 
acts on $X$
by multiplication on the $n$ columns, and we divide by the trivial action of 
the diagonal one-dimensional subtorus. (More precisely, $(z_1,z_2,\dots,z_n) \in (\CC^*)^n$ 
acts on a vector $t^je_i$ by multiplying it by $z_i$\@.  This action on $X$ induces an 
action on $\Gr$\@.)  Note that the lattices corresponding 
to coweights are exactly the torus-fixed points.

This torus action gives rise to a moment map $\momm$ from $\Gr$ to a real
$(n-1)$-dimensional Euclidean space $\vs$\@.  It maps compact irreducible
torus-invariant subvarieties of $\Gr$ onto convex polytopes. It is easy to
describe this map explicitly in terms of lattices:  We first define a
Hermitian inner product on our ambient vector space $X$ by declaring
$\{t^je_i\}$ to be an orthonormal basis.  We also take $\vs$ to be the
hyperspace perpendicular to the vector $(1,1,\dots,1)$ in $\reals^n$ (with
respect to the standard inner product). If $Y=\underline{\lambda}$ is a
torus-fixed point, we may define $\momm(Y)$ to be the orthogonal
projection of $\lambda\in \reals^n$ onto $\vs$\@.  On an arbitrary lattice
$Y$, first write it as an orthogonal sum $Y=\underline{\lambda}\oplus V$
where~$V$ is a finite-dimensional vector space. (For instance, take
$\underline{\lambda}=\cspan\{ t^je_i \mid t^je_i \in Y \}$, and take $V$
to be the orthogonal complement in $Y$\@.)  Choose an orthonormal basis
$v_1,v_2,\dots,v_k$ for $V$\@.  Define $x\in\reals^n$ by $x_i=\sum_j |{\rm
Proj}_i v_j|^2$ where ${\rm Proj}_i v_j$ is the orthogonal projection of
$v_j$ onto the $i\th$ column. 
Then~$\momm(Y)$ is the orthogonal
projection of $\lambda+x \in \reals^n$ onto $\vs$\@.

\subsection{Definition of MV-cycles} 
\label{sec:Definition of MV-cycles}

We recall (an essentially equivalent version of) the Mirkovi\'c-Vilonen definition of MV-cycles.   
Given any lattice $Y$, let us define its {\em highest coweight} $\lambda$ and its {\em lowest 
coweight} $\mu$ to be such that $\underline{\lambda}$ results from using the torus action to 
``flow left'' and $\underline{\mu}$ results from using it to ``flow right''; that is, we 
require that $\lim_{z\rightarrow 0} (1,z,z^2,\dots,z^{n-1})\cdot Y = \underline{\lambda}$ and 
$\lim_{z\rightarrow \infty} (1,z,z^2,\dots,z^{n-1})\cdot Y = \underline{\mu}$\@.  For example, 
in Figure~\ref{fig:lattice-examples}, the first lattice corresponds to the highest coweight 
for the second lattice. Let $S_{\lambda}$ denote the set of lattices with highest coweight 
$\lambda$ and~$T_{\mu}$ denote the set of lattices with lowest coweight $\mu$\@.  
So $S_\lambda$ and $T_\mu$ are stable and unstable manifolds 
with respect to a Morse function given by a certain component of the moment map $\momm$\@.

We define the {\em MV-cycles} with highest coweight $\lambda$ and lowest coweight $\mu$ to 
be the 
irreducible components of the closure of $S_{\lambda} \cap T_{\mu}$\@. 
This intersection is pure dimensional~\cite{MV2}.  

Notice that a coweight $\nu$ shifts $S_\lambda$ and $T_\mu$ to $S_{\lambda+\nu}$ and 
$T_{\mu+\nu}$ respectively. Hence, the MV-cycles with highest 
coweight~$\lambda$ and lowest 
coweight $\mu$ are isomorphic to the MV-cycles with highest 
coweight~$\lambda+\nu$ and lowest
coweight~$\mu+\nu$\@.  We will say that two MV-cycles are equivalent if 
one is a 
shift of another; it will often be convenient to work with the equivalence classes 
rather than with the cycles themselves.

\subsection{Kostant parameter set} 
\label{sec:kostant-parameter} 

Here we recall the definition of Kostant partitions in type~A,  
which are just formal sums of positive roots.  
But it is convenient to develop an unusual notation for these as 
pictures of loops on a Dynkin diagram, for two reasons: 
(1) With this notation, it is easy to see that each Kostant partition 
may naturally be given the structure of a partially ordered set, 
which will be useful.  (2) These pictures make it easy to visualize 
the relationship between Kostant partitions and lattices in the 
affine Grassmannian.

Recall that the Dynkin 
diagram in type $\rm A_{n-1}$ is $n-1$ dots in a row, one for each simple root $\alpha_i$ (connected by line 
segments, which we will not draw). We denote a positive root by a loop around a sequence of 
consecutive dots in the Dynkin diagram, and we call the number of dots it encloses the length 
of the loop.  The simple roots are the loops of length $1$\@. The other positive roots are loops 
of length $\geq2$; each is the sum of the simple roots corresponding to the dots it encloses. 
A {\em Kostant picture} is a picture of the Dynkin diagram together with a 
finite number of such loops. We draw the loops so that if the dots contained in one loop are a 
proper subset of the dots contained in another, then the one is encircled by the other; if two 
loops contain precisely the same dots, we still draw one encircled by the other. In this way, 
the loops in a Kostant picture are partially ordered by encirclement. We write $L\subset L'$ 
if loop~$L'$ encircles loop~$L$, and $L\subseteq L'$ if either $L=L'$ or $L'$ encircles~$L$\@.
The number of loops in a Kostant picture~$\pp$ is denoted $|\pp|$\@.

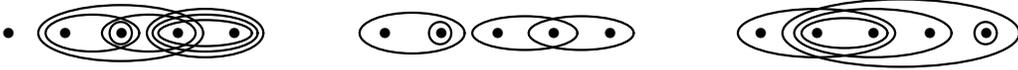
\begin{figure}[ht] 
\begin{center} \begin{picture}(380,35) \setlength{\unitlength}{.5mm}
                                                                                                                 
\def\smalldot{{$\bullet$}}
                                                                                                                 
\put(0,10){\smalldot} \put(15,10){\smalldot} \put(30,10){\smalldot} \put(45,10){\smalldot} 
\put(60,10){\smalldot}
                                                                                                                 
\put(100,10){\smalldot} \put(115,10){\smalldot} \put(160,10){\smalldot} 
\put(130,10){\smalldot} \put(145,10){\smalldot}
                                                                                                                 
\put(200,10){\smalldot} \put(215,10){\smalldot} \put(260,10){\smalldot} 
\put(230,10){\smalldot} \put(245,10){\smalldot}
                                                                                                                 
\thicklines
                                                                                                                 
\put(30,10){\put(1.5,1.6){\ellipse{6}{6}}} \put(22.5,10){\put(1.5,1.6){\ellipse{25}{10}}} 
\put(37.5,10){\put(1.5,1.6){\ellipse{25}{10}}} \put(52.5,10){\put(1.5,1.6){\ellipse{28}{10}}} 
\put(52.5,10){\put(1.5,1.6){\ellipse{25}{7}}} \put(52.5,10){\put(1.5,1.6){\ellipse{31}{13}}} 
\put(30,10){\put(1.5,1.6){\ellipse{44}{15}}}
                                                                                                                 

\put(115,10){\put(1.5,1.6){\ellipse{6}{6}}} \put(152.5,10){\put(1.5,1.6){\ellipse{28}{9}}} 
\put(137.5,10){\put(1.5,1.6){\ellipse{28}{9}}} \put(107.5,10){\put(1.5,1.6){\ellipse{28}{11}}}
                                                                                                                 
                                                                                                                 
\put(260,10){\put(1.5,1.6){\ellipse{6}{6}}} \put(222.5,10){\put(1.5,1.6){\ellipse{23}{8}}} 
\put(237.5,10){\put(1.5,1.6){\ellipse{63}{18}}} \put(230,10){\put(1.5,1.6){\ellipse{42}{13}}} 
\put(215,10){\put(1.5,1.6){\ellipse{42}{13}}}
                                                                                                                 
\end{picture} 
\end{center} 
\caption{Examples of Kostant pictures.} 
\label{fig:Kostant-pictures} 
\end{figure}

The {\em Kostant 
parameter set} $\kost$ is the collection of all Kostant pictures. 
Figure~\ref{fig:Kostant-pictures} shows examples of Kostant pictures 
for $n=6$\@.
If $\alpha_{ij}$ is the root that is the 
sum of simple roots $\alpha_i+\dots+\alpha_{j-1}$, then these three pictures represent the 
following Kostant partitions: $\alpha_{3}+\alpha_{24}+\alpha_{35}+\alpha_{25}+3\alpha_{46}$, 
$\alpha_{2}+\alpha_{13}+\alpha_{35}+\alpha_{46}$ and 
$\alpha_{5}+\alpha_{24}+\alpha_{14}+\alpha_{25}+\alpha_{26}$\@. 

It will be useful to imagine each of the $n-1$ dots of the Dynkin diagram as lying on the 
boundary between two of the $n$ columns of basis vectors for $X$\@.  Then, associated to each 
loop $L$, will be the vector subspace $V_L$ spanned by the columns the loop passes through. To 
be precise, for each $i=1,2,\dots,n$, let $V_i$ be the span of the basis vectors $t^je_i$ in 
the 
$i\th$ column. If we number the dots of the Dynkin diagram $1,2,\dots,n-1$ and if a loop 
$L$ contains dots $\ell,\ell+1,\dots,r-1$ then $V_L=V_{\ell,r}=V_\ell\oplus V_{\ell+1}\oplus 
\dots \oplus V_{r}$\@.  We say that the loop $L$ has its {\em left end} at column~$\ell$, 
its {\em right end} at column~$r$, and {\em passes through} columns $\ell,\ell+1,\dots,r$\@.

By an {\em extended Kostant picture} we mean an ordered pair $\tilde \pp = (\pp, \base)$ where 
$\pp$ is a Kostant picture and $\base$ is a coweight.  (It is sometimes convenient to use this 
terminology: given an extended Kostant picture $\tilde \pp$ and an integer $m$ that's 
less than or equal to 
each component of $\base$, we say that, relative to level $m$, $\tilde \pp$ has $\base_i-m$ 
loops of length zero, or {\em zero loops}, in column $i$\@.) Let $\tilde \kost$ denote the collection of all 
extended Kostant pictures.

Given an extended Kostant picture $\tilde \pp$, we define its {\em highest coweight} $\lambda$ 
and {\em lowest coweight} $\mu$: we set $\lambda_i$ equal to $\base_i$ plus the number of 
loops of $\pp$ with left end at column $i$; we set $\mu_i$ equal to $\base_i$ plus the number 
of loops of $\pp$ with right end at column $i$\@.  Observe that $\pp$ together with any one of 
$\base, \lambda, \mu$ uniquely determines $\tilde \pp$\@.

\subsection{Parametrization of MV-cycles} 
\label{sec:Parametrization of MV-cycles}

We partition $\Gr$ into pieces parametrized by $\tilde \kost$, by defining a map which sends 
$Y\in \Gr$ to $(\pp(Y),\base(Y))\in \tilde\kost$\@. 


Let $\pp(Y)$ be the Kostant picture where the number of loops encircling dots $i,i+1,\dots,j-1$ is
equal to the dimension of the vector space $\big(Y\cap V_{i,j}\big)\big/\big((Y\cap
V_{i+1,j})+(Y\cap V_{i,j-1})\big)$\@. Let $\base(Y)$ be the coweight associated to 
the span of all basis vectors $t^je_i$ contained in $Y$\@.

Then using \cite{AK} terminology we say that the lattice $Y$ 
is {\em weakly compatible} to the extended Kostant picture 
$\tilde \pp = (\pp(Y), \base(Y))$, and we denote by $\MV(\tilde \pp)$ the closure of the set of 
lattices 
weakly compatible to $\tilde \pp$\@.

\begin{theorem} \label{thm:param} {\cite{AK}} The map $\tilde \pp 
\rightarrow \MV(\tilde \pp)$ is a one-to-one correspondence between extended Kostant pictures 
and MV-cycles. \end{theorem}

\noindent
Under the equivalence relation by shifts, we have that 
$\MV((\pp,\base_1))$ is equivalent to $\MV((\pp,\base_2))$ by 
a shift of $\base_2-\base_1$\@.  If we let $\MV(\pp)$ 
denote the equivalence class $\{ \MV(\tilde \pp) \mid \tilde \pp = 
(\pp, \base) \mbox{ for some } \base \}$, then we also have 
a one-to-one correspondence $\pp\rightarrow \MV(\pp)$ between 
Kostant pictures and equivalence classes of MV-cycles.


\begin{remark} Kamnitzer~\cite{Kamnitzer} generalized Theorem~\ref{thm:param} to other 
types. In particular, for any $\gp$ and a choice of a reduced word $\bf i$ for the longest 
element of the Weyl group of $\gp$, he parametrized the set of equivalence classes of MV-cycles 
by the Kostant parameter set.   The parametrization of Theorem~\ref{thm:param} coincides 
with Kamnitzer's parametrization for $\gp=\gln$ and a certain choice of $\bf i$. \end{remark}

The relationship between a Kostant picture and the lattices inside the corresponding MV-cycle 
can perhaps be seen more clearly by choosing a basis for the lattice: Suppose that $\tilde \pp 
= (\pp, \base)$ is an extended Kostant picture and $Y$ is a lattice with basis 
$\{y_L\}_{L\in\pp} \cup \{ t^je_i \mid t^je_i \in \underline{\base} \}$ satisfying $y_L \in 
V_L$ and $t\cdot y_L \in \cspan\{y_{L'} \mid L'\subset L \} + (\underline{\base} \cap V_L)$\@.  
Then it is not hard to see that $Y\in \MV(\tilde \pp)$\@.
Moreover, if $L\in\pp$ has left end at column $\ell$ then the maximum value 
of $j$ for which $y_L$ can have a nonzero component at $t^{-j}e_\ell$ is equal to $\base_\ell + 
|\{ L'\in\pp \mid L'\subseteq L \mbox{ and } L' \mbox{ has left end at column } \ell \}|$;  
similarly for the right end $r$ of $L$\@.  Then $Y$ is weakly compatible to $\tilde \pp$ 
provided that each $y_L$ actually does have a nonzero component at this maximum position, for 
both its left and right end.

For example, let $Y_2$ and $Y_3$ denote the second and third lattices in 
Figure~\ref{fig:lattice-examples} and let $\pp_2$ and $\pp_3$ denote the second and third 
Kostant pictures in Figure~\ref{fig:Kostant-pictures}.  Then $Y_2$ is weakly compatible to 
$\tilde \pp_2 = (\pp_2, (1,-2,-2,-3,-2,-2) )$ and $Y_3$ is weakly compatible to $\tilde \pp_3 
= (\pp_3, (-3,-3,-3,-3,-2,0) )$\@. 

It's also not hard to see from this description that each $\MV(\tilde \pp)$ has dimension 
equal to the sum of the lengths of the loops of $\pp$~\cite{AK}.

\subsection{Richardson varieties} 
\label{sec:richardson varieties}

As an important special case of MV-cycles one finds all of the {\em Richardson varieties} in a 
Grassmannian---i.e. intersections of opposite Schubert varieties.  These are exactly the 
MV-cycles corresponding to Kostant pictures where no loop is encircled by another.

Suppose that $\pp$ is such a Kostant picture; note that its loops are naturally ordered from 
left to right, say $L_1,L_2,\dots, L_{|\pp|}$\@.  Consider the MV-cycle $\MV=\MV(\tilde \pp)$ 
corresponding to $\tilde \pp = (\pp, \base)$, where, to be definite, choose 
$\underline{\base}=X_0$. Let $X_1=t^{-1} \cdot X_0$\@. Then each lattice $Y$ in $\MV$ 
satisfies $X_0 \subseteq Y \subseteq X_1$ with $\cdim (Y/X_0) = |\pp|$\@.  

In fact, the MV-cycle $M$ is the intersection of the two Schubert varieties $S_{\ell}$ and 
$S_{r}$ determined by the left and right ends of the loops of $\pp$\@.  Let $\ell_1 < \ell_2 <  
\dots < \ell_{|\pp|}$ denote the left ends and $r_1 < r_2 < \dots < r_{|\pp|}$ the right ends. 
Then $S_\ell$ is the collection of lattices $Y$ (still with $X_0 \subseteq Y \subseteq X_1$ 
and $\cdim (Y/X_0) = |\pp|$) such that $\cdim (Y\cap V_{\ell_i,n}) / (X_0 \cap V_{\ell_i,n}) 
\geq |\pp|-i+1$  
for $i=1,2,\dots, |\pp|$; and~$S_r$ is the collection of lattices $Y$ 
such that $\cdim (Y\cap V_{1,r_i}) / (X_0 \cap V_{1,r_i}) \geq i$ for $i=1,2,\dots, |\pp|$\@.

\subsection{MV-polytopes} 
\label{sec:MV-polytopes}

MV-polytopes are defined as the moment map images of MV-cycles.  Included among them are all 
of the Weyl polytopes from representation theory (the convex hulls of a Weyl group 
orbit through a coweight).  An MV-polytope $\momm(\MV(\tilde \pp))$ is a translation of a Weyl 
polytope if and only if the number of copies of a particular loop in $\pp$ depends only on the 
length of the loop; for instance if $\pp$ contains five loops of length 1 around the second 
dot of the Dynkin diagram, then it must also contain exactly five loops of length 1 around 
each of the other dots as well.

One reason to be interested in MV-polytopes is that they are closely related to representation 
theory.  For instance, as shown by the first author~\cite{A2}, both weight multiplicities 
and tensor product multiplicities may be 
expressed as the number of MV-polytopes fitting inside a certain region.  Recent work of 
Kamnitzer~\cite{Kamnitzer} shows that this is really an alternative way of looking at 
Bernstein-Zelevinsky's computation of tensor product multiplicities as the number of lattice 
points inside a convex polytope~\cite{BZ}.

In~\cite{AK} all MV-polytopes are computed in type~A, by constructing an explicit map from 
the Weyl group onto the vertices of each polytope.  More recently Kamnitzer~\cite{Kamnitzer} 
has 
generalized this, using Bernstein and Zelevinsky's work~\cite{BZ} to compute MV-polytopes in 
any type.

\section{Convolution of MV-cycles.}
\label{sec:conv}

We will explain how to define a product structure 
on MV-cycles. More specifically, given two MV-cycles $M$ and $N$ we will 
describe a family $M\convC N$ fibered over $\CC$, whose 
fibers over $\cplx^*$ are canonically isomorphic to $M\times N$ and whose fiber over $0$ 
has top-dimensional component equal to  
the union of other MV-cycles $K$ with  multiplicities 
$c^K_{MN}$\@. Then the \convolution product of $M$ and~$N$ is defined to be
\begin{equation}
\label{product}
M\conv N=\sum c_{MN}^K K.
\end{equation}
We will sketch a definition of this product for all types in Section~\ref{sec:BD} and then 
concentrate on a lattice model version of it in type A. Most of this 
section can be omitted on first reading, as we will summarize the 
important results at the beginning of the next section; 
only their statements 
as well as the partial order defined in Section~\ref{sec:order} 
are needed to understand the rest of the paper.

\subsection{\BD Grassmannian and \convolution}
\label{sec:BD}

The key object used in the definition of the \convolution product of
MV-cycles is the \BD Grassmannian $\Gr\convC\Gr$ over~$\cplx$ of the group~$\gp$\@. 
It is an infinite-dimensional 
space which fibers over~$\CC$\@. Just as for the affine Grassmannian, 
the \BD Grassmannian can 
be approximated by 
finite-dimensional algebraic varieties, so that in most cases 
we will be working in a finite-dimensional space.  All 
fibers of $\Gr\convC\Gr$ over $\cplx^*$ are isomorphic to $\Gr\times\Gr$ 
while the fiber over $0$ is
isomorphic to $\Gr$\@. A precise definition of $\Gr\convC\Gr$ along with its
many properties is discussed in \cite{BD}, \cite{MV2}, \cite{Nadler}. 
Instead of repeating this definition we will give a lattice description of 
$\Gr\convC\Gr$ in type A a little later. 

Both the affine Grassmannian and the \BD Grassmannian can be defined for any Borel 
subgroup $\borel$ of $\gp$\@. It turns out that the sets of points of the 
affine Grassmannians for $\gp$ and~$\borel$ are in bijection. However, their topological and 
algebraic structures differ.  In particular, the connected components for the affine 
Grassmannian of $\borel$ are the spaces $S_\lambda$ (similarly $T_\mu$ for the opposite Borel $\borel^-$).
Analogously, the sets of points of the \BD Grassmannians for $\gp$ and $\borel$ are identical, but 
the 
topological and algebraic structures differ. Each irreducible component 
of the \BD Grassmannian for $\borel$ is a family $S_\lambda\convC 
S_\nu$ ($T_\mu\convC T_\gamma$ for~$\borel^-$) which fibers 
over $\CC$ such that its fibers over $\cplx^*$ are naturally isomorphic to 
$S_\lambda\times S_\nu$ (respectively 
$T_\mu\times T_\gamma$) and its fiber over $0$ is isomorphic to 
$S_{\lambda+\nu}$ (respectively $T_{\mu+\gamma}$). We will treat both $S_\lambda\convC
S_\nu$ and $T_\mu\convC T_\gamma$ as subspaces of $\Gr\convC\Gr$\@.

This allows one to define the family $M\convC N$ 
as the irreducible component of 
the closure of the intersection $S_\lambda\convC
S_\nu \cap T_\mu\convC T_\gamma$ 
whose fibers over $\cplx^*$ are $M\times N$\@. Then the fiber of 
$M\convC N$ over~$0$ lies inside the intersection 
$S_{\lambda+\nu}\cap T_{\mu+\gamma}$, 
and its top dimensional components are MV-cycles  for $\overline{S_{\lambda+\nu}\cap 
T_{\mu+\gamma}}$
which are used to 
define the product~(\ref{product}).

If instead of using $M$ and $N$ here 
we used their shifts $\alpha \cdot M$ and $\beta\cdot N$, 
then every fiber of $(\alpha \cdot M)\convC (\beta\cdot N)$ would be 
a shift of the corresponding fiber of $M\convC N$\@. Hence (\ref{product}) defines a 
product on equivalence classes of MV-cycles.

The \convolution product is commutative and associative. 
The commutativity comes from the $\ZZ_2$~symmetry of $\Gr\convC\Gr$ which identifies 
$M\conv_C N$ with $N\convC M$\@.  The associativity 
follows from the fact that it is possible to define the product of 
three MV-cycles using a more general version of the \BD Grassmannian 
fibered over $\CC^2$\@.  (See, for example,  \cite{MV2} or \cite{Nadler}.  
But note that they define 
the \BD Grassmannian in more generality, 
where it is fibered over $C\times C$ for some algebraic curve $C$\@.  
Our $\Gr\convC \Gr$ is a restriction of this 
to a copy of $\CC$ inside $C\times C$ that 
transversely intersects the diagonal of $C\times C$ at zero.)  

\begin{remark}
It is possible to define an action of the torus $\torus$ on $\Gr\convC \Gr$ so that it acts 
diagonally on every fiber over $\CC^*$ and is isomorphic to the usual $\torus$ action at the 
fiber over $0$\@. There is a moment map on $\Gr\convC \Gr$ 
associated to this action.  Then the 
moment map image (and more generally the Duistermaat-Heckman measure) of all fibers of 
$M\convC N$ are the same. Since every fiber over $\CC^*$ is isomorphic to $M\times N$, its 
moment map image is the Minkowski sum of the MV-polytopes 
$\momm(M)$ and $\momm(N)$\@. Hence every MV-polytope $\momm(K)$ for $c_{MN}^K\neq 0$ 
lies within this Minkowski sum. This is a very useful tool in proving that 
certain coefficients from~(\ref{product}) are zero. Also if we use the standard definition 
of convolution of measures on $\mathfrak t_\RR$, the Duistermaat-Heckman measures of 
MV-cycles satisfy an equation exactly like~(\ref{product}).
\end{remark}

\subsection{Lattice model of the \BD Grassmannian in type A} 
For $\gp=\gln$ several lattice models for the \BD 
Grassmannian are known.  One is due to Nadler \cite{Nadler2}.  We 
present another model, which only constructs certain finite-dimensional subfamilies inside of 
$\Gr\convC\Gr$; but since any family $M\convC N$ will embed into one of these, 
it will be sufficient for the study of convolution of MV-cycles. 

Let $\Gr_0$ be the subspace of $\Gr$ containing those lattices $Y$ 
for which $X_0\subset Y$\@. Clearly, each connected component of $\Gr_0$
is a finite-dimensional subvariety of $\Gr$\@.  
We will describe the family $\Gr_0\convC \Gr_0$, which is a subfamily of 
$\Gr\convC\Gr$\@. 
Since any two MV-cycles $M$ and $N$ can be shifted to lie in $\Gr_0$,
their convolution $M\convC N$ can be shifted to lie in $\Gr_0\convC 
\Gr_0$\@. 
Hence the convolution product can be studied within  
$\Gr_0\convC \Gr_0$ instead of the whole \BD Grassmannian.

Let $\widetilde X$ be the vector space of all formal---possibly infinite---linear combinations of 
vectors $t^je_i$ with $-\infty<j<\infty$ and $1\leq i\leq n$\@. As before, 
$X_0$ is the span of all $t^je_i$ with $j\geq 0$, when only finite linear 
combinations are allowed. 
For a complex number $s\neq 0$, we denote the sum
\begin{equation}
\label{power series}
\frac{1}{t}(1+\frac{s}{t}+\frac{s^2}{t^2}+\cdots
+\frac{s^k}{t^k}+\cdots)
\end{equation} 
by $\frac{1}{t-s}$\@.
Let $X^s$ be the span of $X_0$ and the
vectors $\frac{e_i}{(t-s)^{j}}$ for $j> 0$ (where we treat $\frac{1}{(t-s)^{j}}$ as 
the~$j\th$ power of the series~(\ref{power series})). A subspace 
$Y$ of $X^s$  will be called an $s$-lattice
if it is invariant under multiplication by $(t-s)$, it contains $X_0$, and
for a large $N$ it is contained inside the space spanned by $X_0$ and
the vectors $\frac{e_i}{(t-s)^{j}}$  for $N\geq j> 0$\@.

Now consider the variety $\mathcal F$ of all the subspaces $Y$ of
$X$ of the form $Y=\yz+\ys$, where $\yz$ is a $0$-lattice (namely an 
ordinary 
lattice containing $X_0$) and
$\ys$ is an $s$-lattice for $s\neq 0$\@. Clearly, $\mathcal F$ is a union 
of 
finite-dimensional algebraic varieties and an $s$-fibration over 
$\CC-\{0\}$ whose fibers are 
isomorphic to $\Gr_0\times \Gr_0$\@.

Finally, consider the space $\mathcal X$ of subspaces $Z$ of $\tilde X$ 
with $X_0\subset Z$ and $\dim (Z/X_0)<\infty$\@. The space~$\Gr_0\convC\Gr_0$ is now defined 
to be the 
closure of  $\mathcal F$ within $\mathcal X$.

\subsection{Example}
\label{sec:example}
This lattice model for the \BD Grassmannian is useful for finding 
concrete families $M\convC N$ which deform $M\times N$ into unions of 
MV-cycles. Sometimes the lattice model can be used to guess correctly 
which coefficients on the right side of the equation~(\ref{product}) 
are not zero. However, it is not well suited for computing these 
coefficients. In particular it is very hard to see using this model 
whether a nonzero coefficient is equal to $1$ or to some other integer. Let us 
present an illustration.

Assume $n=3$\@. Let $\pp$ be the Kostant partition $\alpha_1$ and $\qq$ be 
the Kostant partition~$\alpha_2$\@. Define $\tilde\pp =(\pp,\underline{0})$, 
$\tilde\qq =(\qq,\underline{0})$, where $\underline{0}=(0,0,0)$\@. Then
$\MV(\tilde\pp)\conv \MV(\tilde\qq)=\MV(\tilde\rr)+\MV(\tilde\ss)$, where 
$\tilde \rr=(\alpha_1+\alpha_2,\underline{0})$ and 
$\tilde \ss=(\alpha_{13},(0,1,0))$\@. In 
terms of the Kostant pictures we can write this equation as
\begin{equation}
\label{fig:example}
\begin{picture}(200,0)

\put(0,10){ \setlength{\unitlength}{.5mm}}

\def\smalldot{{\tiny $\bullet$}}

\put(20,0){\smalldot} \put(30,0){\smalldot}

\put(20,0){\put(1.7,1.8){\ellipse{6}{6}}}

\put(5,-1){$M($}
\put(35,-1){$)\conv M($}

\put(60,0){\smalldot} \put(70,0){\smalldot}

\put(70,0){\put(1.7,1.8){\ellipse{6}{6}}}

\put(75,-1){$)=M($}

\put(110,0){\smalldot} \put(120,0){\smalldot}

\put(110,0){\put(1.7,1.8){\ellipse{6}{6}}}

\put(120,0){\put(1.7,1.8){\ellipse{6}{6}}}

\put(125,-1){$)+M($}

\put(160,0){\smalldot} \put(170,0){\smalldot}

\put(160,0){\put(6.2,1.6){\ellipse{17}{8}}}

\put(175,-1){$).$}

\end{picture}
\end{equation}

Both MV-cycles $\MV(\tilde \pp)$ and $\MV(\tilde \qq)$ are projective lines 
$\CC\PP^1$, while $\MV(\tilde \rr)$ and $\MV(\tilde \ss)$ are projective 
planes $\CC\PP^2$ intersecting along a $\CC\PP^1$\@. Hence the family 
$M\convC N$ degenerates $\CC\PP^1\times \CC\PP^1$ to two copies of 
$\CC\PP^2$ intersecting along $\CC\PP^1$\@. There are many 
ways of visualizing this family. For example, let $W$ be a four-dimensional 
vector space with basis $w_1,w_2,w_3,w_4$\@. Then the $s$-fiber of this 
family can be identified with the pair of lines, one sitting inside the 
span of $w_1$ and $w_2$ and the other inside the span of $w_2+sw_4$ and 
$w_3$\@.  The zero fiber contains those two-dimensional planes 
inside $W$ 
which either contain the vector $w_2$ or sit inside the span of 
$w_1,w_2,w_3$\@.  

We will use the lattice model to show that both 
MV-cycles appearing on the right side of 
equation~(\ref{fig:example}) should indeed be there. However, checking 
that both multiplicities are $1$ and that there are no other 
terms on the right side requires other 
techniques (for example, 
moment map images or, more generally, Duistermaat-Heckman measures of 
MV-cycles). 

Consider the space $\mathcal F(\tilde \pp,\tilde\qq)$ which contains all 
the vector subspaces 
$Y\subset \widetilde X$ spanned by $X_0$ and two vectors $y_1$ and $y_2$ 
of the form $y_1=at^{-1}e_1+bt^{-1}e_2$ and 
$y_2=\frac{ce_2}{t-s}+\frac{de_3}{t-s}$ 
for some $s\neq 0$\@. The closure of $\mathcal F(\tilde \pp,\tilde\qq)$ is 
the family 
$\MV(\tilde \pp)\convC  \MV(\tilde \qq)$\@. 

To show that a certain MV-cycle $K$ enters the product $\MV(\tilde \pp)\conv 
\MV(\tilde \qq)$ 
with nontrivial coefficient it is enough to construct a one-parameter family 
$Y(s)$ inside  $\MV(\tilde \pp)\convC \MV(\tilde \qq)$ such 
that $Y(s)$ lies inside the $s$-fiber of $\MV(\tilde \pp)\convC \MV(\tilde 
\qq)$, and $Y(0)$ 
is a generic point of $K$\@. We will construct two such families in our 
example, one for each MV-cycle on the right side of 
equation~(\ref{fig:example}).

Given any numbers $a,b,c,d$, we can define the family $Y(s)$ to be the  
span of $X_0$, $y_1(s)=at^{-1}e_1+bt^{-1}e_2$ 
and $y_2(s)=\frac{ce_2}{t-s}+\frac{de_3}{t-s}$\@. 
Then $Y(0)$ is the span of 
$X_0$, $y_1(0)=at^{-1}e_1+bt^{-1}e_2$, and $y_2(0)=ct^{-1}e_2+dt^{-1}e_3$, 
which is a generic point in $\MV(\tilde \rr)$\@.

Similarly, given numbers $a,d$ and $b\neq 0$, let $Y(s)$ be spanned by 
$X_0$, $y_1(s)=sat^{-1}e_1+bt^{-1}e_2$, and 
$y_2(s)=\frac{be_2}{t-s}+s\frac{de_3}{t-s}$ for $s\neq 0$\@. Then $Y(0)$ is 
the 
span of $X_0$, $y_1(0)=bt^{-1}e_2$, and 
$\lim_{s\to 0} \frac{y_1(s)}{s}-\frac{y_2(s)}{s}= 
at^{-1}e_1-bt^{-2}e_2-dt^{-1}e_3$, which is a generic point in 
$\MV(\tilde\ss)$\@.

\subsection{A partial order on MV-cycles and its relation to the \convolution 
product.}
\label{sec:order}
In this section we present a property of the convolution of MV-cycles. Our 
statement of this property is specific to type A and the proof
uses the lattice model for the \BD 
Grassmannian.  However, it should be possible to use the 
techniques developed by Kamnitzer~\cite{Kamnitzer} 
to generalize this to other types, as will be discussed below.

Let us start by defining a partial order on extended Kostant pictures. 
For an extended Kostant picture $\tilde \pp =(\pp, \base)$, fix an $m$ 
such that each $\base_i>m$, and think of $\tilde \pp$ as having 
$\base_i-m$ zero loops in column~$i$\@.
Consider loops $A,B$ in $\tilde \pp$ which overlap, 
in the sense that neither encircles the other and they  
pass through at least one common column. 
Define a new Kostant picture $\tilde \pp'$ by replacing $A$ and $B$ by 
their union $A\cup B$ (i.e.\ the loop passing though the columns which 
either $A$ or $B$ pass through) and intersection $A\cap B$
(i.e.\ the loop passing though the columns which both $A$ and $B$
pass through; if there is only one such column, say $i$, then this is a zero 
loop, in which case $\base_i$ is increased by 1). 
Call this operation $\tilde\pp \rightarrow \tilde\pp'$ the {\it fusion} of loops $A$ 
and~$B$\@. 
Notice that fusion does 
not change the highest and lowest coweights of the extended Kostant picture.

For example, consider the extended Kostant picture $\tilde \rr$ from 
Section~\ref{sec:example}.  
It contains two length one loops and
their fusion is 
a length two loop and a length zero loop.  So this fusion 
produces the extended Kostant picture $\tilde \ss$\@.

Define a partial order on $\tilde \KK$ by saying that 
$\tilde\pp\geq \tilde \qq$ if $\tilde \qq$ is produced from $\tilde\pp$ 
by a sequence of fusions. 
To see that this partial order is well-defined it is enough to show that 
a sequence of nontrivial fusions of an extended Kostant picture $\tilde \pp$ cannot 
produce 
$\tilde \pp$ again.
Fix an integer $m$ as before.
For $1\leq i\leq j\leq n$ define~$d_{ij}(\tilde \pp)$ to be the number of
loops (including the zero loops) of $\tilde\pp$ which lie between columns~$i$ and~$j$; in other 
words $d_{ij}(\tilde \pp)$ is the number of loops of $\tilde \pp$ which do 
not pass through any of the columns $1,\dots, i-1,j+1,\dots, n$\@. 
Let $D(\tilde \pp)$ denote the collection of numbers $d_{ij}(\tilde \pp)$\@. We say 
that $D(\tilde \pp)\leq D(\tilde \qq)$ if $d_{ij}(\tilde\pp)\leq 
d_{ij}(\tilde\qq)$ for every~$i,j$\@; and $D(\tilde \pp)< D(\tilde \qq)$ if 
$D(\tilde \pp)\leq D(\tilde \qq)$ and 
$d_{ij}(\tilde\pp)<
d_{ij}(\tilde\qq)$ for some~$i,j$\@. 

It is clear that if $\tilde \qq$ is produced out of $\tilde \pp$ by a 
sequence of nontrivial fusions, then $D(\tilde\qq)> D(\tilde\pp)$, which 
proves that the partial order on $\tilde \KK$ is well-defined. 
Surprisingly the converse is also true.

\begin{lemma}  
\label{lem:order}
Assume $\tilde\pp$ and $\tilde \qq$ have the 
same highest and lowest coweights.   Fix a common $m$ to define the 
numbers $d_{ij}(\tilde \pp)$ and $d_{ij}(\tilde \qq)$\@. Then $D(\tilde\qq)> 
D(\tilde\pp)$  if and only if $\tilde\qq<\tilde\pp$\@.
\end{lemma}

\begin{proof} It remains to show that $D(\tilde\qq)> D(\tilde\pp)$  
implies $\tilde\qq<\tilde\pp$\@. If $\tilde\qq$ and $\tilde\pp$ have a 
common loop, we can remove it and proceed by induction on the number of 
loops. 
If all the loops of $\tilde \qq$ and $\tilde\pp$ are different, pick a
loop $C$ inside $\tilde \qq$ which does not encircle any other loop. 
Let the left end of $C$ be $\ell$ and the right end $r$\@. 
Consider all the loops of $\tilde \pp$ which encircle $C$, and 
among these, consider only those that do not encircle any other 
loops which encircle $C$; these are ordered in a sequence from left to right.  
Moreover there must be at least two of them: in particular, since the lowest and
highest coweights of $\tilde \pp$ and $\tilde \qq$ are the same, 
there is one with right end $r$ and another with left end $\ell$.
Let $A$ and $B$ be any two consecutive loops in this sequence, 
and let $\tilde \pp'$ be produced from $\tilde \pp$ by the fusion of $A$ and $B$\@.  
Then it is not hard to check that
$D(\tilde\qq)\geq D(\tilde\pp')$\@.  If $D(\tilde\qq)= D(\tilde\pp')$ 
then $\tilde \qq=\tilde \pp'$ and we are done; otherwise $D(\tilde\qq) > D(\tilde\pp')$
and we proceed by induction with $\tilde \pp$ replaced by $\tilde \pp'$\@.
\end{proof}


Let $\pp$ and $\qq$ be two Kostant pictures. Define $\pp+\qq$ to be the 
Kostant picture containing the loops of both $\pp$ and $\qq$; in other 
words, we just add the Kostant partitions.  
If $\tilde \pp=(\pp,\base)$ and $\tilde 
\qq=(\qq,\base')$, then define $\tilde \pp+\tilde 
\qq=(\pp+\qq,\base+\base')$\@.

\begin{prop}
\label{prop:order}
If $\MV(\tilde \pp)\conv \MV(\tilde \qq)=\sum c_{\tilde \rr}\MV(\tilde \rr)$ 
then 
$c_{\tilde \rr}=0$ for  $\tilde\rr > \tilde\pp+\tilde\qq$\@. 
\end{prop}

\begin{proof}
Without loss of generality we may assume that $\tilde \pp=(\pp,\underline 
0)$ and $\tilde \qq=(\qq,\underline 0)$\@. For the rest of the proof fix 
$m=0$ when defining the $d_{ij}$'s for extended Kostant pictures.

For a lattice $Y$ containing $X_0$ define $D(Y)$ to be the set of numbers 
$d_{ij}(Y)=\dim(Y\cap V_{ij}/X_0\cap V_{ij})$\@. (Recall that 
$V_{ij}$ is the span of columns $i$ through $j$ inside $X$\@.)  
It is easy to see that $D(Y)\geq D(\tilde \pp)$ whenever $Y\in \MV(\tilde 
\pp)$ and $D(Y)\geq D(\tilde \qq)$ whenever $Y\in \MV(\tilde \qq)$\@. Moreover, 
generically these inequalities are equalities.

Now let $Y$ be a point in $\MV(\tilde \pp)\convC \MV(\tilde \qq)$\@. We can 
again define $D(Y)$ to be the set of numbers $d_{ij}(Y)=\dim(Y\cap 
\widetilde V_{ij}/X_0\cap \widetilde V_{ij})$, where $\widetilde V_{ij}$ 
is now 
the span of columns 
$i$ though $j$ of  $\widetilde X$\@. If $Y$ is in a fiber of 
$M\convC N$ over $\cplx^*$ then $Y=\yz+\ys$ with $\yz\in \MV(\tilde \pp)$ 
and $\ys\in 
\MV(\tilde \qq)$\@.  Hence
$$
D(Y)=D(\yz)+D(\ys)\geq 
D(\tilde\pp)+D(\tilde\qq)=
D(\tilde\pp+\tilde\qq).$$
Since the functions $d_{ij}$ are lower semi-continuous, we conclude that 
$D(Y)\geq D(\tilde\pp+\tilde\qq)$ 
for any $Y\in \MV(\tilde \pp)\convC \MV(\tilde 
\qq)$, and in particular for $Y$ in the zero fiber. Applying 
Lemma~\ref{lem:order} finishes the proof.
\end{proof}

\begin{remark} We believe that Proposition~\ref{prop:order} holds in a more 
general situation. Namely as shown by Kamnitzer~\cite{Kamnitzer}, once a 
reduced word for the longest element of the Weyl group of~$\gp$ is chosen, 
it is possible to parametrize MV-cycles by extended Kostant partitions. 
Moreover, 
Kamnitzer constructed functions which naturally generalize the numbers 
$d_{ij}$ in the above proof. So, we expect that 
Proposition~\ref{prop:order} holds for any 
type as soon as the proper partial order on extended Kostant partitions is defined. 
The proof of the more general statement should  rely on the 
semi-continuity of Kamnitzer's analogues of the numbers $d_{ij}$\@. 
\end{remark}

\subsection{The leading term of the convolution product}
The converse of Proposition~\ref{prop:order} is false: 
one may have $c_{\tilde \rr} = 0$ even when 
$\tilde\rr < \tilde\pp+\tilde\qq$\@.   
But not for 
$\tilde\rr = \tilde\pp+\tilde\qq$:

\begin{prop}
\label{prop:lead}
If $\MV(\tilde \pp)\conv \MV(\tilde \qq)=\sum c_{\tilde \rr}\MV(\tilde \rr)$ 
then $c_{\tilde\pp+\tilde\qq}\neq 0$\@.
\end{prop}

\begin{conj}
\label{conj:lead}
The coefficient $c_{\tilde\pp+\tilde\qq}$ is always equal to $1$\@.
\end{conj}

\noindent
The rest of this section is devoted to the proof of 
Proposition~\ref{prop:lead} using the lattice model. 

Without loss of generality 
$\tilde \pp=(\pp,\underline 0)$ and $\tilde \qq=(\qq,\underline 0)$\@. 
Let $Z$ be a lattice weakly compatible to 
$\tilde \pp+\tilde \qq$\@. It will be enough to prove that for each 
such $Z$ there exists a one-parameter family $Z(s)$ inside 
$\MV(\tilde\pp)\convC \MV(\tilde \qq)$ with $Z(0)=Z$ and $Z(s)$ sitting inside 
the $s$-fiber of $\MV(\tilde\pp)\convC \MV(\tilde \qq)$\@.

The lattice $Z$ has a basis $\{y_L\}_{L\in\tilde \pp+\tilde \qq}$  
indexed by the loops of $\tilde \pp+\tilde \qq$ as described in Section~\ref{sec:Parametrization of 
MV-cycles}. 
(Here we treat each vector $t^je_i$ as the basis vector associated to a zero loop.)
This basis satisfies the property that each vector $ty_L$ lies in the span 
of the vectors associated to loops encircled by $L$\@. So we can aways 
uniquely write
\begin{equation}
\label{eq:relation}
ty_L=\sum_{L'\subset L} a_{L}^{L'}y_{L'}.
\end{equation}
Notice that if we count all the zero loops there are infinitely many loops 
encircled by $L$; nevertheless the summation on the right side 
of~(\ref{eq:relation}) is finite.

Let us now decompose the Kostant picture 
$\pp+\qq$ into $\pp$ and $\qq$: color every loop in 
$\pp+\qq$ either red or green, 
so that the red loops give $\pp$ and the green loops give $\qq$\@.

Define vectors $z_L=z_L(s)$ inside $\widetilde X$ by
$$
z_L=\frac{1}{t-s} \sum_{L'\subset L} a_{L}^{L'}z_{L'}
$$
if $L$ is red and by 
$$
z_L=\frac{1}{t} \sum_{L'\subset L} a_{L}^{L'}z_{L'}
$$
if $L$ is green.

Let $Z(s)$ be the span of all the vectors $z_L(s)$ and $X_0$\@. Clearly, 
$Z=Z(0)$\@. So, to prove Proposition~\ref{prop:lead} it is enough to show 
that $Z(s)$ is inside the $s$-fiber of $\MV(\tilde\pp)\convC \MV(\tilde \qq)$ 
for $s\neq 0$\@.

To show this we will construct vectors $x_L=x_L(s)$ of the form
\begin{equation}
\label{eq:x}
x_L=z_L+\sum_{L'\subset L} b_L^{L'}z_{L'}
\end{equation}
such that for every red loop $L$, the vector $(t-s)x_L$ is a linear 
combination of vectors~$x_{L'}$ 
for $L'\subset L$ a red loop or a zero loop.  
Similarly, for every green 
loop~$L$, the vector 
$tx_L$ will be a linear
combination of vectors $x_{L'}$
for $L'\subset L$ a green loop or a zero loop. 
Then we define~$\ys(s)$ to be spanned by red vectors $x_L$ 
and $X_0$, and 
$\yz(s)$ to be spanned 
by green vectors~$x_L$ and $X_0$\@. 
Clearly, $Z(s)=\yz(s)+\ys(s)$, $\ys(s)$ 
is 
an $s$-lattice inside $\MV(\tilde \pp)$, and $\yz(s)$ is a $0$-lattice 
inside 
$\MV(\tilde \qq)$, proving that $Z(s)$ is inside the $s$-fiber of 
$\MV(\tilde\pp)\convC \MV(\tilde \qq)$\@.  It remains to construct vectors~$x_L$ 
with the above properties.

Construct the vectors $x_L$ by induction on inclusion of 
loops. Namely, assume we can construct such vectors for every loop 
$L'\subset L$ and we will prove that $x_L$ exists for $L$\@. 
Assume that $L$ is red; an analogous argument can be given if $L$ is green.

Start by setting $x=z_L$ and write
\begin{equation}
\label{eq:xx}
(t-s)x=\sum_{L'\subset L} c^{L'}x_{L'}.
\end{equation}
If in the above equation $c^{L'}=0$ for all green loops $L'$ we can just 
set $x_L=x$\@. Otherwise, let $L^g$ be one of the largest green loops with 
$c^{L^g}\neq 0$; in other words if $L^g$ is encircled by a green loop $L'$ 
then $c^{L'}=0$\@. We can then modify $x$ to guarantee $c^{L^g}=0$ and 
proceed by induction. Namely we replace $x$ by the vector
$$
x'= x +\frac{c^{L_g}}{s} x_{L^g}.
$$
Since $L^g\subset L$, the vector $x'$ is still of the form~(\ref{eq:x}). 
And we have
\begin{align*}
(t-s) x'&= (t-s)x +\frac{(t-s)c^{L^g}}{s} x_{L^g}
=(t-s)x + \frac{tc^{L^g}}{s} x_{L^g} 
- c^{L^g}x_{L^g} \\
&=\sum_{L'\subset L} c^{L'}x_{L'} -c^{L^g}x_{L^g} 
+\frac{c^{L^g}}{s} tx_{L^g}
\end{align*}
Thus the vector $x_{L^g}$, as well as all green vectors $x_{L'}$ for loops 
$L'$ encircling $L^g$, have coefficient zero when we write $(t-s)x'$ as a 
linear combination of the $x_{L'}$\@. This allows us to proceed by induction and 
construct the vector $x$ of the form~(\ref{eq:x}) such that in equation~(\ref{eq:xx}),
$c^{L'}=0$ whenever $L'$ is green. Then we can set $x_L=x$\@.

\section{Convolution algebra and $\palg$}
\label{sec:convolution algebra}

\subsection{Properties of the convolution product} Let us recall the results of 
Section~\ref{sec:conv} as well as restate them in terms of Kostant pictures 
rather than extended Kostant pictures. We will assume that $\gp=\sln$ or 
$\gp=\gln$ for most of this section.

Equation~(\ref{product}) is invariant under shifts. Hence we can define 
the convolution product on equivalence classes of MV-cycles
\begin{equation}
\label{eq:product}
\MV(\pp)\conv \MV(\qq)=\sum c_{\pp,\qq}^\rr \MV(\rr).
\end{equation}
From this point on, 
whenever we mention MV-cycles we will really mean 
their equivalence classes, 
and we will only use the convolution 
defined in equation~(\ref{eq:product}).
Also, the partial order on extended Kostant pictures immediately 
induces one on Kostant pictures: 
we say $\pp>\qq$ if $\tilde\pp=(\pp,\base) > \tilde\qq=(\qq,\base')$ for some $\base, \base'$\@.

Summarizing Section~\ref{sec:conv} (in particular 
Propositions~\ref{prop:order} and \ref{prop:lead}) 
we have the following properties of convolution:
\begin{enumerate}
\item {The convolution product is commutative and associative,}
\item 
\label{property2} {$c_{\pp,\qq}^\rr=0$ for $\rr > \pp+\qq$,}
\item 
\label{property3} {$c_{\pp,\qq}^{\pp+\qq}\neq 0$.}
\end{enumerate}

\noindent
The commutativity and associativity (which hold for any $\gp$)   
imply
that the $\CC$-span of MV-cycles has an algebra structure.  
We call this the convolution algebra $\calg$ for the group $\gp$\@.

\subsection{The isomorphism.} For a group $\gp$, denote by $\unip$ the unipotent 
radical of a Borel subgroup, and let $\palg$ be the algebra of 
functions on $\unip$\@. 
Assume $\gp$ is either $\sln$ or $\gln$, so that $\unip$ can be 
identified with upper triangular $n\times n$ matrices with $1$'s along the 
diagonal; then $\palg$ is the polynomial ring $\CC[\xx]$ 
in the matrix entries 
$\xx=\{ x_{ij} \}_{1\leq i<j\leq n}$
above the diagonal.

Let us define a homomorphism $\isom : \palg\to \calg$ by sending $x_{ij}$ to 
$\MV(\pp_{ij})$, where $\pp_{ij}$ is the Kostant picture with a single 
loop passing through columns $i$ through $j$.

\begin{theorem}
\label{thm}
The map $\isom$ is an isomorphism.
\end{theorem}

\begin{proof}
For a loop $L$ passing through columns $i$ through $j$ let $x_L=x_{ij}$, 
and for a Kostant partition~$\pp$ define the monomial $x^\pp=\prod_{L\in \pp} 
x_L$\@. Clearly the monomials $x^\pp$ provide a basis for $\palg$\@.

Let us show that $\isom$ is surjective. We will prove by induction that 
every 
$\MV(\pp)$ 
is in the image of~$\isom$\@.  Assume  $\MV(\qq)\in {\rm Im}(\isom)$  for 
every $\qq < \pp$, and we'll prove the same thing for $\pp$\@. By 
properties~(\ref{property2}) and~(\ref{property3}) of the convolution product
$$
\isom(x^\pp)=\prod_{L\in \pp}\isom (x_L)= a \MV(\pp) + \sum_{\qq<\pp} 
a_\qq \MV(\qq)
$$
with $a\neq 0$\@. Using the induction assumption, we conclude that $\MV(\pp)\in 
{\rm Im}(\isom)$, which proves that $\isom$ is surjective.

Let us show that $\isom$ is injective. For an integer $k$, let 
$\palg^k$ be the span of the monomials $x^\pp$ with~$|\pp|\leq k$; 
similarly, let 
$\calg^k$ be the span of the $\MV(\pp)$ with $|\pp|\leq k$\@. Clearly $\dim 
(\palg^k)=\dim (\calg^k)$\@. Moreover, by property~(\ref{property2}) of the 
\convolution product, $\isom^{-1}(\calg^k)\subset \palg^k$\@. Hence $\isom$ is 
injective.
\end{proof}

Let us denote $\isom^{-1}(\MV(\pp))$ by $\poly_\pp$ and call it an 
MV-polynomial. 
It is easy to see that 
$$
\poly_\pp=\sum_{\qq\leq \pp} c_{\qq}^{\pp}x^\qq
$$
with $c_\pp^\pp\neq 0$\@.
Computation of the the coefficients $c_{\qq}^\pp$ is a challenging 
unsolved problem, which we hope has a combinatorial solution. 
Section~\ref{sec:determinant} is devoted to stating some conjectures about this.

The MV-polynomials form a basis of $\palg$, which is expected 
to coincide with one of the known canonical bases of $\palg$\@. We are not 
experts in this field and will not speculate which one it is.

\begin{remarks} 

$\bullet$  The first author in~\cite{Anderson-thesis} defined a conjectural Hopf algebra structure 
on $\calg$ and conjectured that $\palg$ and $\calg$ are isomorphic as Hopf algebras.
Theorem~\ref{thm} proves part of that conjecture in type~A. 

$\bullet$ Feigin, Finkelberg, Kuznetzov, and Mirkovi\'c~\cite{FM, FFKM} gave a geometric 
construction of the universal enveloping algebra dual to $\palg$. 

$\bullet$  We believe that Theorem~\ref{thm} can be generalized to 
other types by using Kamnitzer's results~\cite{Kamnitzer}. As mentioned 
before, we expect that analogues of properties~(\ref{property2}) 
and~(\ref{property3}) of the convolution product for other 
types can be proved using Kamnitzer's results. Then the above proof of 
Theorem~\ref{thm} should  work for other types with only minor 
modifications.

$\bullet$ The partial order on Kostant partitions, 
viewed as a partial order on the monomials of $\CC[\xx]$, 
is a diagonal term order.

$\bullet$ Property~(\ref{property2}) of \convolution is very 
similar to the multiplicative property of  Lusztig's dual canonical basis proved 
by Caldero~\cite{Caldero}. Since MV-polynomials are expected to form one of the known 
canonical basis of $\palg$, this similarity is not surprising. 
\end{remarks}

\section{Determinantal formulas for MV-polynomials}
\label{sec:determinant}

A large unsolved problem 
is to find an explicit formula for every MV-polynomial.  
This is equivalent to the problem of finding a formula for the convolution 
of any two MV-cycles (at least insofar 
as a solution to one problem gives an inductive method for solving the other).  
Likely, these are difficult problems; all we will 
be able to do here is to give a mysterious 
conjectural formula for some MV-polynomials as certain determinants. Even more mysterious is the 
connection to the Fomin-Zelevinsky theory of cluster algebras, for our conjectures 
apply only to those MV-polynomials which are also cluster monomials. 

\subsection{MV-polynomials for Richardson varieties}

Recall that according to the isomorphism of Section~\ref{sec:convolution algebra}, 
if $\pp$ is a Kostant picture consisting of a single loop $L$ encircling 
dots $i,i+1,\dots,j-1$ of the Dynkin diagram, then the corresponding 
MV-polynomial is just the matrix entry $x_{ij}$\@.

Then, by computing convolutions of MV-cycles, one may calculate MV-polynomials 
for more complicated Kostant pictures.  
For instance, from the computation of the example in~\ref{sec:example}, we see that 
\begin{equation*}
\begin{picture}(200,0)
                                                                                
\put(0,10){ \setlength{\unitlength}{.5mm}}
                                                                                
\def\smalldot{{\tiny $\bullet$}}

\put(10,-5){\smalldot} \put(18,-5){\smalldot}
                                                                                
\put(10,-5){\put(1.7,1.8){\ellipse{6}{6}}}
                                                                                
\put(18,-5){\put(1.7,1.8){\ellipse{6}{6}}}
                                                                                
\put(6,-1){$\poly$}
\put(30,-1){$ \! =  \poly$}
                                                                                
\put(45,-4){\smalldot} \put(53,-4){\smalldot}
                                                                                
\put(45,-4){\put(1.7,1.8){\ellipse{6}{6}}}
                                                                                
\put(60,-1){$\poly$}

\put(65,-4){\smalldot} \put(73,-4){\smalldot}
                                                                                
\put(73,-4){\put(1.7,1.8){\ellipse{6}{6}}}

\put(80,-1){$- \!\ \poly$}
                                                                                
\put(97,-4){\smalldot} \put(105,-4){\smalldot}
                                                                                
\put(97,-4){\put(5.2,1.6){\ellipse{15}{8}}}
                                                                                
\put(111,-1){$=x_{12}x_{23}-x_{13}$}
                                                                                
\end{picture}
\end{equation*}
\vspace{0.5mm}

By doing such computations, one quickly finds 
that the MV-polynomial for any of the Richardson varieties 
is given by a minor.  More precisely, suppose $\pp$ is a Kostant picture 
without any loops encircled by other loops.  
Let $\ell_1,\ell_2,\dots,\ell_{|\pp|}$ and 
$r_1,r_2,\dots,r_{|\pp|}$ denote the left and right ends 
of the loops of $\pp$, ordered from left to right.  
Then $\poly_\pp$ is the minor $\det ( x_{\ell_i,r_j} )$\@.  
Here $(x_{ij})\in \unip$ denotes an 
arbitrary element of the group of upper triangular unipotent matrices, 
so that $x_{ii}=1$ and $x_{ij}=0$ if $i>j$. 
Moreover, all the nonvanishing minors arise in this way.
For instance, in the above example, we have 
\begin{equation*}
\begin{picture}(200,0)
                                                                                
\put(0,10){ \setlength{\unitlength}{.5mm}}
                                                                                
\def\smalldot{{\tiny $\bullet$}}

\put(12,-4){\smalldot} \put(20,-4){\smalldot}
                                                                                
\put(12,-4){\put(1.7,1.8){\ellipse{6}{6}}}
                                                                                
\put(20,-4){\put(1.7,1.8){\ellipse{6}{6}}}
                                                                                
\put(6,-1){$\poly$}
\put(30,-1){$=$}
\put(40,-1){$\left|
\begin{array}{cc}
x_{12} & x_{13}  \\
1      & x_{23}  \\
\end{array}
\right|$}
\end{picture}
\end{equation*}
\vspace{0.5mm}

\subsection{Which MV-polynomials are determinants?}
\label{sec:whicharedeterminants}

As one inductively computes more and more complicated MV-polynomials, 
one observes that often they are naturally expressible as determinants; 
indeed it is difficult to find one that is not.  We begin by 
saying precisely what we mean by ``naturally expressible''.

Suppose $\poly_\pp$ is an MV-polynomial, with leading term $x^\pp$\@.  
Label the loops of $\pp$ in any order: $L_1,L_2,\dots,L_{|\pp|}$\@.  
We will construct a $|\pp|$ by $|\pp|$ matrix~$A_\pp$ whose 
rows and columns are labelled by the loops $L_1,L_2,\dots,L_{|\pp|}$ 
and the product of whose diagonal entries is $x^\pp$\@.  
We set the $ij\th$ entry of~$A_\pp$ equal to $x_{\ell_i,r_j}$ where 
$\ell_i$ is the left end of $L_i$ and $r_j$ is the right end of $L_j$ 
(recalling that $x_{ii}=1$ and $x_{ij}=0$ if $i>j$).

One observes that for many Kostant pictures $\pp$ 
the MV-polynomial $\poly_\pp$ is equal to the determinant 
of a matrix $\hat A_\pp$, obtained from $A_\pp$ by changing some of its 
entries to zeros.  In particular, this is true 
for the Richardson varieties, in which case no entries are 
changed and $\hat A_\pp=A_\pp$, the determinant of which is 
the minor in the previous section.

\vspace{5mm}
\noindent
{\bf Example.} Let $\pp$ be the Kostant picture 

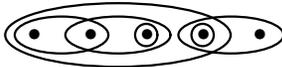
\begin{picture}(380,35) 
\setlength{\unitlength}{.5mm}
                                                                                                                 
\def\smalldot{{$\bullet$}}
                                                                                                                 
\put(0,10){\smalldot} \put(15,10){\smalldot} \put(30,10){\smalldot} \put(45,10){\smalldot} 
\put(60,10){\smalldot}
                                                                                                                 
\thicklines

\put(30,10){\put(1.5,1.6){\ellipse{6}{6}}}                                                                                                                 
\put(45,10){\put(1.5,1.6){\ellipse{6}{6}}} 
\put(7.5,10){\put(1.5,1.6){\ellipse{25}{10}}} 
\put(22.5,10){\put(1.5,1.6){\ellipse{28}{10}}} 
\put(52.5,10){\put(1.5,1.6){\ellipse{28}{10}}} 
\put(22.5,10){\put(1.5,1.6){\ellipse{60}{17}}}

\end{picture} 

\noindent
One can compute that the corresponding MV-polynomial is 
\begin{eqnarray*}
\poly_\pp &=& x_{13}x_{24}x_{34}x_{15}x_{45}x_{46}
-x_{13}x_{24}x_{34}x_{45}^2x_{16}
-x_{13}x_{24}x_{15}x_{35}x_{46}
+x_{13}x_{24}x_{35}x_{45}x_{16} \\
& &-x_{13}x_{34}x_{15}x_{45}x_{26}
+x_{13}x_{34}x_{25}x_{45}x_{16}
+x_{13}x_{15}x_{35}x_{26}
-x_{13}x_{25}x_{35}x_{16} \\
& &-x_{23}x_{14}x_{34}x_{15}x_{45}x_{46}
+x_{23}x_{14}x_{34}x_{45}^2x_{16}
+x_{23}x_{14}x_{15}x_{35}x_{46}
-x_{23}x_{14}x_{35}x_{45}x_{16} \\ 
& &-x_{14}x_{15}x_{25}x_{46}
+x_{14}x_{15}x_{45}x_{26}
+x_{24}x_{15}^2x_{46} 
-x_{24}x_{15}x_{45}x_{16} 
-x_{15}^2x_{26}
+x_{15}x_{25}x_{16}
\end{eqnarray*}

Let us order the six loops of $\pp$ 
according to the locations of their left ends 
in the drawing of the Kostant picture, 
working from left to right (so $L_1$ is the loop of length $4$ 
and $L_6$ is the rightmost loop of length $1$).  
Having chosen this order, we may write down the matrix $A_\pp$, 
and (for instance by trial and error) we find that 
$\poly_\pp=\det(\hat A_\pp)$ after changing ten entries to zero.
\[
A_{\pp}
=
\left(
\begin{array}{cccccc}
x_{15} & x_{13} & x_{14} & x_{14} & x_{16} & x_{15} \\
x_{15} & x_{13} & x_{14} & x_{14} & x_{16} & x_{15} \\
x_{25} & x_{23} & x_{24} & x_{24} & x_{26} & x_{25} \\
x_{35} &      1 & x_{34} & x_{34} & x_{36} & x_{35} \\
x_{45} &      0 &      1 &      1 & x_{46} & x_{45} \\
x_{45} &      0 &      1 &      1 & x_{46} & x_{45} \\
\end{array}
\right) \hspace{10mm}
\hat A_{\pp}=\left(
\begin{array}{cccccc}
x_{15} & \bz    & \bz    & \bz    & x_{16} & \bz    \\
x_{15} & x_{13} & x_{14} & x_{14} & x_{16} & x_{15} \\
x_{25} & x_{23} & x_{24} & x_{24} & x_{26} & x_{25} \\
\bz    &      1 & \bz    & x_{34} & \bz    & x_{35} \\
x_{45} &      0 &      1 &      1 & x_{46} & x_{45} \\
\bz    &      0 & \bz    &      1 & \bz    & x_{45} \\
\end{array}
\right)
\]

Notice that there might be several different ways to 
set entries equal to zero so as to get the same determinant.  
For instance, in this example, subtracting the third 
column from the fourth column in $\hat A_\pp$ 
would result in setting an additional three entries equal 
to zero without changing the determinant.  
Nevertheless, as we will discuss in the next section, 
we believe that in general there should be a 
canonical best choice of how to do this.


\vspace{5mm}
\noindent
{\bf Counterexample.} 
It is not easy to find an example of an MV-polynomial 
that cannot be expressed as a determinant in the way just described.  
The first counterexample is in type $\rm A_5$, and is actually the only 
such example we have been able to compute.  
The Kostant picture $\pp_{!}$ is 

\begin{picture}(380,35) 
\setlength{\unitlength}{.5mm}
                                                                                                                 
\def\smalldot{{$\bullet$}}
                                                                                                                 
\put(0,10){\smalldot} \put(15,10){\smalldot} \put(30,10){\smalldot} \put(45,10){\smalldot} 
\put(60,10){\smalldot}
                                                                                                                 
\thicklines

\put(30,10){\put(1.5,1.6){\ellipse{6}{6}}}                                                                                                                 
\put(7.5,10){\put(1.5,1.6){\ellipse{25}{10}}} 
\put(52.5,10){\put(1.5,1.6){\ellipse{25}{10}}} 
\put(30,10){\put(1.5,1.6){\ellipse{44}{13}}} 

\end{picture} 

\noindent
and the MV-polynomial is 
\begin{eqnarray*}
\poly_{\pp_{!}} &=& x_{34}x_{13}x_{46}x_{25}
-x_{23}x_{34}x_{46}x_{15}
-x_{34}x_{45}x_{13}x_{26}
+x_{23}x_{34}x_{45}x_{16} 
-x_{13}x_{25}x_{36}
-x_{46}x_{14}x_{25} \\
& &+x_{24}x_{46}x_{15}
+x_{13}x_{35}x_{26} 
+x_{23}x_{36}x_{15}
+x_{45}x_{14}x_{26}
-x_{23}x_{35}x_{16}
-x_{45}x_{24}x_{16} \\
& &-2x_{15}x_{26}
+2x_{25}x_{16}
\end{eqnarray*}

\noindent
Then, ordering the loops left to right by their left ends, 
\[
A_{\pp_!}=\left(
\begin{array}{cccc}
x_{13} & x_{15} & x_{14} & x_{16} \\
x_{23} & x_{25} & x_{24} & x_{26} \\
1      & x_{35} & x_{34} & x_{36} \\
0      & x_{45} &      1 & x_{46} \\
\end{array}
\right)
\]

\noindent
and one can verify that $\poly_{\pp_!}$ is not 
equal to the determinant of any $\hat A_{\pp_!}$ 
obtained by changing some entries to zeros.

We believe there to be a relation between the Fomin-Zelevinsky theory of 
cluster algebras~\cite{FZ1,FZ2} and the question of which MV-polynomials are 
expressible as determinants.  
The algebra $\palg$ carries 
the structure of a cluster algebra.  
One of the basic ingredients of this structure is a canonical 
set of generators for $\palg$, called {\em cluster variables}, and
a superset of {\em cluster monomials}, which are certain monomials in the 
cluster variables.  (Of course a cluster monomial is not necessarily a 
{\em monomial} in the variables $\xx$; it is a polynomial in these variables.)
For completeness, more details are provided in Appendix~\ref{appendix}; 
in particular, we will give a precise definition of the cluster monomials 
in $\palg$, and describe the algorithm we used to compute them.  
In this section, 
we simply view cluster algebra theory as a black box that outputs 
certain polynomials---the cluster monomials---which we conjecture 
are a subset of MV-polynomials.

\begin{conj}
\label{conj:cluster} 
Every cluster monomial is an MV-polynomial that is naturally expressible as a determinant;   
that is, if $\psi$ is a cluster monomial, then there 
is a Kostant picture $\pp$ such that $\psi=\poly_\pp={\rm det}(\hat A_\pp)$, 
for some matrix $\hat A_\pp$ obtained from $A_\pp$ by setting 
some entries equal to zero.
\end{conj}

\begin{remark}
If this conjecture is true then the above counterexample $\poly_{\pp_!}$ 
is an MV-polynomial which is not a cluster monomial.
\end{remark}

\noindent
This conjecture rests on three pieces of evidence:  

\begin{enumerate}

\item For every cluster monomial that was small enough so as 
to be able to also compute the MV-polynomial with the same leading term, 
these two polynomials were identical.  

\item We have been able to find the matrix $\hat A_\pp$ 
for approximately a hundred different cluster monomials in type $\rm A_5$\@.  
A Kostant picture with nine or ten loops is about the maximum 
size for which it is usually possible to find $\hat A_\pp$ within a few hours 
using educated guesswork (and letting a computer do the algebra of course).  
A real problem is that we do not know an efficient algorithm for finding $\hat A_\pp$\@.

\item We strongly believe that the above counterexample of an MV-polynomial 
that is not expressible as a determinant is not a cluster monomial, 
because it has not appeared in the first 719 cluster variables for type $\rm A_5$ 
output by the computer~\cite{AK*}, whereas all other Kostant pictures of this size 
are accounted for.   This same polynomial was used by Leclerc~\cite{Leclerc2} to give a 
counterexample to a conjecture of Berenstein-Zelevinsky. Leclerc showed that this 
MV-polynomial is an element of the dual canonical basis, but its square is not. 

\end{enumerate}

\begin{remark}
Although Conjecture~\ref{conj:cluster} only refers to the cluster monomials, 
one might hope that a more general formula, of the form
$$
\poly_\pp=\sum_{w\in {\rm Perm}(\pp)} {\rm sgn}(w) c_w x^w 
$$
holds for any MV-polynomial $\poly_\pp$\@.  Here ${\rm Perm}(\pp)$ is the group 
of permutations of the loops of $\pp$, ${\rm sgn}(w)$ is the sign of the permutation, 
and $x^w=\prod_{L\in \pp} x_{\ell_{w(L)},r_L}$ where $w(L)$ has left end at $\ell_{w(L)}$
and $L$ has right end at $r_L$; as always $x_{ii}=1$ and $x_{ij}=0$ if $i<j$\@.  
(Note that different $w$ can give the same $x^w$\@.)  
The coefficients $c_w$ would be nonnegative integers somehow determined  
combinatorially from $\pp$ and the permutations $w$\@.  
(Unfortunately we are far from understanding the $c_w$---all we can really 
say is that $c_w$ should be $0$ if $x^w=x^\qq$ with $\qq>\pp$\@.)  
In the case that $\poly_\pp$ is a cluster monomial, 
each $c_w$ would be $0$ or $1$, and the sum would reduce to the 
definition of a determinant $\det(\hat A_\pp)$\@.
\end{remark}

\subsection{Which entries should be changed to 0?}

Here, we will give a 
conjecture as to how to find~$\hat A_\pp$ for a moderately large 
class of Kostant pictures $\pp$ for which we believe $\poly_\pp$ to be a 
determinant; that is, we will specify exactly which entries of $A_\pp$ are 
to be changed to zeros.

First, let us describe the class of Kostant pictures we will consider.  
Given any Kostant picture $\pp$, we construct a directed graph $G_\pp$ 
whose vertices are the loops of $\pp$\@.  It will be directed in the 
sense that for each edge, we will say that the loop at one end is on the 
left and the loop at the other end is on the right.

The first step in constructing $G_\pp$ is to draw a ``line diagram''
of $\pp$.  
This will be somewhat similar to the pictures back in Figure~\ref{fig:lattice-examples}. 
For each loop $L$ of $\pp$ and each column $i$ through which it passes, 
let us say that its height $h_i$ is equal to the length 
of the longest chain of loops of $\pp$ ending in $L$, 
say $L_1 \subset L_2 \subset \dots \subset L_h=L$\@.  
Then, starting with a horizontal base of dots at height zero, 
we draw, for each loop in $\pp$, a dot at height $h_i$ for each column 
it passes through, and connect these dots by straight line segments.  

For instance, for the Kostant picture in the first example of this section, 

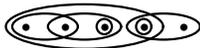
\begin{picture}(300,35) 

\put(35,5){
\setlength{\unitlength}{.35mm}
                                                                                                                 
\def\smalldot{{\tiny$\bullet$}}
                                                                                                                 
\put(0,10){\smalldot} \put(15,10){\smalldot} \put(30,10){\smalldot} \put(45,10){\smalldot} 
\put(60,10){\smalldot}
                                                                                                                 
\thicklines

\put(30,10){\put(1.5,1.6){\ellipse{6}{6}}}                                                                                                                 
\put(45,10){\put(1.5,1.6){\ellipse{6}{6}}} 
\put(7.5,10){\put(1.5,1.6){\ellipse{25}{10}}} 
\put(22.5,10){\put(1.5,1.6){\ellipse{28}{10}}} 
\put(52.5,10){\put(1.5,1.6){\ellipse{28}{10}}} 
\put(22.5,10){\put(1.5,1.6){\ellipse{60}{17}}}

}
\end{picture}

\noindent
we draw

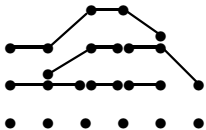
\begin{picture}(380,50) 
\setlength{\unitlength}{.5mm}

\def\smalldot{{$\bullet$}} \def\tinydot{{\tiny$\bullet$}}
                                                                                
\thicklines
                                                                                
\put(20,0){\smalldot} \put(30,0){\smalldot} \put(40,0){\smalldot} \put(50,0){\smalldot} 
\put(60,0){\smalldot} \put(70,0){\smalldot}

\put(20,10){\smalldot} \put(30,10){\smalldot} \put(38.5,10){\smalldot} \put(22,12){\line(1,0){18.5}}

\put(41.5,10){\smalldot} \put(48.5,10){\smalldot} \put(43.5,12){\line(1,0){7}}

\put(51.5,10){\smalldot} \put(60,10){\smalldot} \put(53.5,12){\line(1,0){8.5}}

\put(30,13){\smalldot} \put(41.5,20){\smalldot} \put(32,15){\line(23,14){10}} 
\put(48.5,20){\smalldot} \put(43.5,22){\line(1,0){7}}

\put(51.5,20){\smalldot} \put(60,20){\smalldot} \put(53.5,22){\line(1,0){8.5}} 
\put(70,10){\smalldot} \put(62,22){\line(1,-1){10}}

\put(20,20){\smalldot} \put(30,20){\smalldot} \put(22,22){\line(1,0){10}} 
\put(41.5,30){\smalldot} \put(32,22){\line(10,9){10}}
\put(50,30){\smalldot} \put(43.5,32){\line(1,0){8.5}}
\put(60,23){\smalldot} \put(52,32){\line(10,-7){10}}





\end{picture} 

\noindent (Note that it is easy to draw these line diagrams inductively:  (1) start 
by drawing the line segments for the loops that don't contain other loops; (2) 
for any loop $L$ that encircles only loops for which the line segments have 
already been drawn, place the dot in each column at the lowest position that's 
one unit higher than any dot already drawn in that column {\em for a loop that's 
encircled by} $L$\@.)

Now, let us suppose $L_1$ and $L_2$ are overlapping loops of $\pp$,  
meaning that neither loop encircles the other, 
and they pass through at least one column in common; without loss of 
generality suppose that their left and right ends satisfy 
$\ell_1<\ell_2\leq r_1<r_2$\@.  
We write $L_1\diredge L_2$ provided that 
the line segments we drew for $L_1$ intersect those we drew for $L_2$; 
more precisely, we require that 
in column~$r_1$, the height of~$L_1$ is 
less than or equal to the height of~$L_2$, 
and in column~$\ell_2$, the height of~$L_2$ is 
less than or equal to the height of $L_1$\@.  
Then we connect~$L_1$ and~$L_2$ by an edge in~$G_\pp$, 
with~$L_1$ on the left and~$L_2$ on the right,  
provided that in addition there 
is no third loop~$L$ in~$\pp$ with $L_1\diredge L\diredge L_2$\@.

In the above example, the graph looks like this:

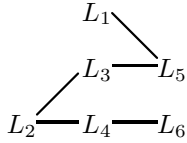
\begin{picture}(380,60) 
\setlength{\unitlength}{.5mm}
                                                                                
\thicklines
                                                                                
\put(25,0){$L_2$} 
\put(45,0){$L_4$} 
\put(65,0){$L_6$} 
\put(45,15){$L_3$} 
\put(65,15){$L_5$} 
\put(45,30){$L_1$} 

\put(33,3){\line(1,0){11}} 
\put(53,3){\line(1,0){12}} 
\put(53,18){\line(1,0){12}} 
\put(33,5){\line(1,1){11}} 
\put(53,32){\line(1,-1){12}} 

\end{picture} 

\noindent
(The labelling of the loops is the same as it was before; 
working left to right, the loops of length 1 are $L_4,L_6$, 
the loops of length 2 are $L_2,L_3,L_5$, and the loop of length 4 is $L_1$.)   
The picture is drawn such that for any edge the name of the loop 
at the left end is situated on the page to the left of the name of the loop 
at the right end.  

Now we may say which Kostant pictures $\pp$ we will consider: 
those for which the graph $G_\pp$ is acyclic, i.e., is 
a disjoint union of trees.  (Notice that the graph for the 
counterexample $\pp_!$ is not one of these; it is a 4-cycle.)
For any such $\pp$ we will conjecture a matrix $\hat A_\pp$ 
for which $\poly_\pp=\det(\hat A_\pp)$\@.  So we must specify 
which entries of $A_\pp$ to change to zeros.  To do this, 
we need to consider paths in the graph $G_\pp$\@.

Suppose that $(L_{a_1},L_{a_2},\dots,L_{a_k})$ is a path of distinct vertices in $G_\pp$ 
(where we are allowed to travel in either direction along an edge).  
We say that it is an {\em allowable} path provided that 
for any three consecutive vertices $L_{a_i},L_{a_{i+1}},L_{a_{i+2}}$ on it: 
\begin{enumerate}
\item If $L_{a_i}\diredge L_{a_{i+1}}$ and $L_{a_{i+2}}\diredge L_{a_{i+1}}$ then $L_{a_i}$ encircles 
$L_{a_{i+2}}$.  
\item If $L_{a_{i+1}}\diredge L_{a_i}$ and $L_{a_{i+1}}\diredge L_{a_{i+2}}$ then $L_{a_{i+2}}$ encircles 
$L_{a_i}$.  
\end{enumerate}

\begin{conj}
\label{conj:zeros} 
Suppose $\pp$ is a Kostant picture for which the graph $G_\pp$ is acyclic.  
Order the loops of $\pp$ and define the matrix $A_\pp$ as before.  
Define $\hat A_\pp$ to be the matrix obtained from $A_\pp$ by changing 
the $ij\th$ entry to $0$ if there is no allowable path in $G_\pp$ from the 
$j\th$ to the $i\th$ loop of $\pp$\@. Then $\poly_\pp=\det(\hat A_\pp)$\@.
\end{conj}

One can check in the above example that this results in changing to zeros 
the ten entries we saw earlier.  For instance, the entry of $\hat A_\pp$ in 
row 1 and column 2 is set to $\mathbf 0$ since there is no allowable path 
from $L_2$ to $L_1$; indeed the path $(L_2,L_3,L_5,L_1)$ is not allowable 
since condition~(1) fails for the final three vertices $(L_3,L_5,L_1)$\@.

\begin{remarks} 

$\bullet$ This conjecture is true for many other Kostant pictures as well, 
but for a different graph, which we don't know how to define in general.  
For instance, if $\pp$ is the Kostant picture

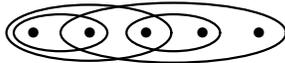
\begin{picture}(380,35) 
\setlength{\unitlength}{.5mm}
                                                                                                                 
\def\smalldot{{$\bullet$}}
                                                                                                                 
\put(0,10){\smalldot} \put(15,10){\smalldot} \put(30,10){\smalldot} \put(45,10){\smalldot} 
\put(60,10){\smalldot}
                                                                                                                 
\thicklines

\put(7.5,10){\put(1.5,1.6){\ellipse{25}{10}}} 
\put(37.5,10){\put(1.5,1.6){\ellipse{25}{10}}} 
\put(15,10){\put(1.5,1.6){\ellipse{44}{15}}} 
\put(37.5,10){\put(1.5,1.6){\ellipse{60}{16}}} 


\end{picture} 

\noindent
then the conjecture does not apply since $G_\pp$ is a 4-cycle.  
Nevertheless, we still have that $\poly_\pp=\det(\hat A_\pp)$, 
if instead of using 
$G_\pp$ to construct $\hat A_\pp$, we use a different graph: 
the one obtained from $G_\pp$ 
by removing the edge joining the left loop of length 2 
to the loop of length~4.

$\bullet$ We believe that the allowable path condition, although correct, is not 
the natural thing to say; 
that is, there should be a different, more general, condition which happens 
to give exactly this answer for this particular class of Kostant pictures. 
\end{remarks}

\subsection{Connecting MV-polynomials to geometry}

We will state a conjecture that directly relates each MV-polynomial 
to the corresponding MV-cycle, in the case where it is a cluster monomial.  
Our motivation was to provide a non-inductive geometric 
definition of each MV-polynomial, but we still cannot do this.

When told this conjecture, a geometer typically 
expresses some combination of perplexity and dismay, 
perhaps asking ``Where does this come from?'' and ``What does it mean?''.  
The answers are: It is a purely empirical conjecture based on examples, 
and, unfortunately, we have no idea what it really means; we are not 
hiding any geometric intuition from the reader.

To state the conjecture, it is convenient to use the natural group action of 
$\gln(\polytinv)$ on the affine Grassmanian $\Gr$:  $\gln$ acts 
by the 
standard representation on $e_1,e_2,\dots,e_n$ and $t$ acts by sending 
each basis vector $t^je_i$ to $t^{j+1}e_i$; extend this by linearity.  
Let $\unip$ and $\unip^-$ denote the subgroups of~$\gln$ 
of upper and lower triangular unipotent matrices.  
Then it is easy to see that in the definition of MV-cycles 
as irreducible components of $\overline{S_\lambda \cap T_\mu}$ 
in Section~\ref{sec:Definition of MV-cycles}, 
we have $S_\lambda=\unip^-(\polytinv)\cdot\underline{\lambda}$ 
and~$T_\mu=\unip(\polytinv)\cdot\underline{\mu}$\@.  
So, given an MV-cycle $M$, with lowest coweight $\mu$, 
there is a dense subset~$M^\circ$ of points $Y\in M$ for which 
we can write $Y=\nel \cdot \underline{\mu}$, where $\nel \in \unip(\polytinv)$; 
of course $\nel$ is not uniquely determined by $Y$\@.

Let us suppose $M=\MV(\tilde \pp)$ where $\tilde \pp = (\pp, \base)$; 
as usual, $\lambda$ and $\mu$ denote the highest and lowest coweights.    
We will define two functions $\fone, \ftwo: M^\circ \rightarrow \cplx$ 
and conjecture that they are the same function.  (Actually, the first function 
$\fone$ will only be well-defined when the MV-polynomial $\poly_\pp$ is a 
cluster monomial, 
and even this well-definedness will be conjectural.)  
The definition of $\fone$ will depend on $\poly_\pp$ 
but the definition of $\ftwo$ will not.  

\subsubsection{The first function}

Let $Y$ be an arbitrary point of $M^\circ$ and write $Y=\nel \cdot \underline \mu$\@.  
In short, we will define $\fone(Y)$ to be the lowest non-zero term 
of $\poly_\pp(\nel)\in \polytinv$, the evaluation of $\poly_\pp$ on 
$\nel\in\unip(\polytinv)$\@.  
Of course $\poly_\pp(\nel)$ depends not just on $Y\in M^\circ$, 
but on $\nel$ as well.  Let $d$ be the integer that is the net 
difference in heights $j$ of the basis vectors $t^{-j} e_i$ of $\underline{\lambda}$ 
and those of $\underline{\mu}$, that is, 
$d=(\sum_{t^{-j}e_i \in \underline{\mu} \mbox{ \tiny but } t^{-j}e_i \notin \underline{\lambda}} \, j) -
(\sum_{t^{-j}e_i \in \underline{\lambda} \mbox{ \tiny but } t^{-j}e_i \notin \underline{\mu}} \, j)$\@.  
Then in the case where  $\poly_\pp$ is a cluster monomial we conjecture 
that the coefficient of $t^d$ in $\poly_\pp(\nel)$ does not depend on the choice of $\nel$, 
and we denote this coefficient by $\fone(Y)$\@.  
It is important to assume that $\poly_\pp$ is a cluster monomial; 
in particular, for the counterexample Kostant picture $\pp_!$, 
this coefficient does depend on the choice of~$\nel$\@. 
As for the other coefficients of $\poly_\pp(\nel)$, 
for $k<d$ we conjecture that the coefficient of $t^k$ is $0$; 
and for $k>d$, any example will show that
the coefficient of $t^k$ depends on the choice of $\nel$\@.


\subsubsection{The second function}

Here $\pp$ is allowed to be any Kostant picture, and we again 
let $Y$ be an arbitrary element of $M^\circ$\@.  
Let $\pi_\lambda: Y \rightarrow \underline{\lambda}$ and 
$\pi_\mu: Y \rightarrow \underline{\mu}$ 
denote orthogonal projection.  
It is not hard to see that $Y=\nel \cdot \underline{\mu}$ 
implies that $\pi_\mu$ is invertible, 
so that we can define 
$\pi = \pi_\lambda \circ \pi_\mu^{-1}:\underline{\mu}\rightarrow\underline{\lambda}$\@.  
Although formally this is a linear map between infinite-dimensional 
vector spaces, it is really a map between finite-dimensional vector spaces of dimension $|\pp|$: 
From the discussion in Section~\ref{sec:Parametrization of MV-cycles}
of bases for lattices inside an MV-cycle, we can write 
$\underline{\lambda}=\underline{\base} \oplus \underline{\hat \lambda}$ and 
$\underline{\mu}=\underline{\base} \oplus \underline{\hat \mu}$, 
where $\underline{\hat \lambda}$ and $\underline{\hat \mu}$ 
are $|\pp|$-dimensional subspaces of $\underline{\lambda}$ and $\underline{\mu}$ 
orthogonal to $\underline{\base}$, 
and $\pi$ restricts to a map $\hat \pi: \underline{\hat\mu} \rightarrow \underline{\hat\lambda}$\@.  
We define~$\ftwo(Y)$ to be the Jacobian of this linear map.  Ordinarily, 
the Jacobian of a map between different vector spaces is only defined up to sign; 
but here the sign is determined since the Kostant picture~$\pp$ 
gives a canonical identification between the basis vectors of 
$\underline{\hat\lambda}$ and of $\underline{\hat\mu}$: 
those for $\underline{\hat\lambda}$ correspond to the left ends of the loops and
those for $\underline{\hat\mu}$ correspond to the right ends of the loops.

\begin{conj}
\label{conj:samefunction} 
Suppose $\pp$ is a Kostant picture for which the MV-polynomial $\poly_\pp$ 
is a cluster monomial.  
Let $M=\MV(\tilde \pp)$ be the corresponding MV-cycle with lowest coweight $\mu$, and let 
$M^\circ$ be its dense subset $M \cap (\unip(\polytinv)\cdot 
\underline{\mu})$\@.  
Then $\fone$ is well-defined and the functions $\fone,\ftwo: M^\circ \rightarrow \cplx$ are identical.
\end{conj}

\begin{remarks}

$\bullet$ It should be possible to extend $\fone=\ftwo$ from $M^\circ$ to a function 
$\chi: M \rightarrow \CC\PP^1$\@. 

$\bullet$ There should be a direct relation between Conjecture~\ref{conj:samefunction} 
and Conjecture~\ref{conj:zeros}.  In particular, one might be able to 
see from the definitions of $\fone$ and $\ftwo$  
exactly which entries of the matrix $A_\pp$ 
must be changed to zeros, because $\fone$ is related to $\det (\hat A_\pp)$ 
and $\ftwo$ is related (via lattices and linear algebra) to the combinatorics of $\pp$\@.  
But we have not been able to see the connection, 
and the conjectures were arrived at independently.

$\bullet$ The evidence for Conjecture~\ref{conj:samefunction} is weaker than that   
for Conjecture~\ref{conj:zeros}; it has not been tested on Kostant 
pictures with more than four or five loops. 
\end{remarks}

\vspace{5mm}
\noindent
{\bf Example.} 
We will verify the conjecture when the Kostant picture $\pp$ is 

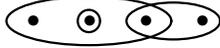
\begin{picture}(380,35) 
\setlength{\unitlength}{.5mm}
                                                                                                                 
\def\smalldot{{$\bullet$}}
                                                                                                                 
\put(0,10){\smalldot} \put(15,10){\smalldot} \put(30,10){\smalldot} \put(45,10){\smalldot}
                                                                                                                 
\thicklines

\put(15,10){\put(1.5,1.6){\ellipse{6}{6}}}                                                                                                                 
\put(37.5,10){\put(1.5,1.6){\ellipse{25}{10}}} 
\put(15,10){\put(1.5,1.6){\ellipse{44}{13}}} 

\end{picture} 

\noindent
In this case, the MV-polynomial is 
\[
\poly_\pp=\left|
\begin{array}{ccc}
x_{14} & \bz    & x_{15} \\
x_{24} & x_{23} & x_{25} \\
x_{34} & 1      & x_{35} \\
\end{array}
\right|
=
x_{23}x_{35}x_{14}-x_{23}x_{34}x_{15}-x_{14}x_{25}+x_{24}x_{15}
\]

Let us take $M=\MV(\pp,\base_0)$ where $\underline {\base_0} = X_0$ 
is spanned by all $t^je_i$ for $j\geq 0$\@.
Then $\underline{\mu}=\underline{\base_0} \oplus \cspan \{ t^{-1}e_3, t^{-1}e_4, t^{-1}e_5 \}$\@.
Write 
\[
\nel=\left(
\begin{array}{ccccc}
1          & a           & e           & h           & j           \\
0          & 1           & b           & f           & i           \\
0          & 0           & 1           & c           & g           \\
0          & 0           & 0           & 1           & d           \\
0          & 0           & 0           & 0           & 1           \\
\end{array}
\right)
\]
for $\nel \in \unip(\polytinv)$, 
where $a = \sum_{k=-\infty}^{\infty} a_k t^k, 
b = \sum_{k=-\infty}^{\infty} b_k t^k, \dots \in \polytinv$\@.
Using the description in Section~\ref{sec:Definition of MV-cycles} of lattices in an MV-cycle, 
it is not hard to verify 
that necessary and sufficient conditions for $Y=\nel \cdot \underline{\mu}$ 
to be in $M^\circ$ are: 
(1)~$e$ has degree $\geq 1$; 
$a,b,d,h,j$ have degrees $\geq 0$;
$c,f,g,i$ have degrees $\geq -1$ 
(where the degree is the smallest $k$ for which the coefficient of $t^k$ is nonzero).  
(2)~The following minors vanish:
\[
\left|
\begin{array}{cc}
b_{0}  & f_{-1} \\
1      & c_{-1} \\
\end{array}
\right|
\hspace{7mm}
\left|
\begin{array}{cc}
b_{0}  & i_{-1} \\
1      & g_{-1} \\
\end{array}
\right|
\hspace{7mm}
\left|
\begin{array}{cc}
h_{0}  & j_{0} \\
f_{-1}  & i_{-1} \\
\end{array}
\right|
\hspace{7mm}
\left|
\begin{array}{cc}
h_{0}  & j_{0} \\
c_{-1}  & g_{-1} \\
\end{array}
\right|
\hspace{7mm}
\left|
\begin{array}{cc}
f_{-1}  & i_{-1} \\
c_{-1}  & g_{-1} \\
\end{array}
\right|
\]
To compute $\fone(Y)$ we calculate 
\begin{align*}
\poly_\pp(\nel) =&\ bgh-bcj-hi+fj \\
=&  
\left(
b_0 
\left|
\begin{array}{cc}
h_{0}  & j_{0} \\
c_{-1} & g_{-1} \\
\end{array}
\right|
-
\left|
\begin{array}{cc}
h_{0}  & j_{0} \\
f_{-1} & i_{-1} \\
\end{array}
\right|
\right) t^{-1} \\
  & 
+
\left(
b_0 g_0 h_0 - b_0 c_0 j_0 - h_0 i_0 + f_0 j_0 
+
b_1
\left|
\begin{array}{cc}
h_{0}  & j_{0} \\
c_{-1}  & g_{-1} \\
\end{array}
\right|
-j_1
\left|
\begin{array}{cc}
b_{0}  & f_{-1} \\
1      & c_{-1} \\
\end{array}
\right|
+h_1
\left|
\begin{array}{cc}
b_{0}  & i_{-1} \\
1      & g_{-1} \\
\end{array}
\right|
\right) \\
 & + \mbox{  {\em terms of degree $1$ and higher}}
\end{align*}
where we have grouped together certain pairs of terms into vanishing minors; therefore
$$\fone(Y)= b_0 g_0 h_0 - b_0 c_0 j_0 - h_0 i_0 + f_0 j_0$$

Now we compute $\ftwo(Y)$\@.  We have 
$\hat\lambda = \cspan \{ t^{-1}e_1, t^{-1}e_2, t^{-1}e_3 \}$,
$\hat\mu = \cspan \{ t^{-1}e_4, t^{-1}e_3, t^{-1}e_5 \}$, 
where we view these as ordered bases in the order written 
(so that the $i\th$~vector corresponds to the $i\th$ loop when 
the loops are ordered left to right by their left ends).  
With respect to these ordered bases, it is easy to check that we have 
\[
\hat\pi=
\left(
\begin{array}{ccc}
h_0 & 0   & j_0\\
f_0 & b_0 & i_0 \\
c_0 & 1   & g_0 \\
\end{array}
\right)
\left(
\begin{array}{ccc}
1   & 0   & d_0 \\
c_0 & 1   & g_0 \\
0   & 0   & 1   \\
\end{array}
\right)^{\mbox{-1}}
\]
The Jacobian of $\hat\pi$ is then the quotient of the determinants of these two matrices, which is 
$$\ftwo(Y)= b_0 g_0 h_0 - b_0 c_0 j_0 - h_0 i_0 + f_0 j_0$$
and we have verified that $\fone(Y)=\ftwo(Y)$ for any $Y\in M^\circ$\@.

\appendix\section{Cluster Algebra}
\label{appendix}

Here we define a very interesting set of generators for $\palg$ 
giving it the structure of a cluster algebra.  
Cluster algebras were discovered by Fomin and Zelevinsky \cite{FZ1,
FZ2, BFZ}, and have proved to be related to many different topics. By now the
theory of cluster algebras is so large, 
that it would be pointless to even try to summarize their motivation and 
many applications here.  All we will do is  
recall the definitions and results needed to define a cluster
algebra inside $\palg$ for $\gp=\gln$ or $\sln$\@. For the more
general treatment of these results see~\cite{BFZ}.

For us a cluster $\cluster$ will be a collection
$(c_1,\dots,c_{\frac{n(n-1)}{2}}, B)$ of
$\frac{n(n-1)}{2}$ elements $c_k$ of $\palg$ and an
$\frac{n(n-1)}{2}\times \frac{(n-1)(n-2)}{2}$ matrix of integers $B$\@.  
The clusters will be defined inductively: 
we will define an original cluster $\cluster_0$ and specify a
mutation rule
designed to produce more clusters from existing clusters. More
specifically,
for a cluster $\cluster$ and an element $c_k$ of this cluster with $k\leq
\frac{(n-1)(n-2)}{2}$ we will
define another cluster
$\mu_k(\cluster)=(c_1,\dots,c'_k,\dots,c_{\frac{n(n-1)}{2}}, \mu_k(B))$
by replacing $c_k$ by a new
variable~$c'_k$ and replacing $B$ by another matrix $\mu_k(B)$\@. The
mutation rule will
be defined in such a way that $\mu_k(\mu_k(\cluster))=\cluster$\@.
Then every cluster will be produced out of the original cluster by a
sequence
of mutations.
                                                                               
The original cluster is defined as follows. Since we can identify $\unip$ with
the space of upper triangular unipotent matrices, every
minor of such an $n\times n$ matrix can be though of as an element of $\palg$\@.
For $1\leq a<b\leq n$ define
$\Delta_{ab}\in \palg$ to be the minor with columns
$a,\dots, b-1$
and the top $b-a$ rows of the $n\times n$ matrix. The functions of the original
cluster are given by  $\Delta_{ab}$\@. More specifically, let
$k(a,b)=\frac{b(b-1)}{2}-a$, and set
$c_{k(a,b)}=\Delta_{ab}$\@. The matrix $B$ is given as follows.
The entry of $B$ in the intersection of row  $k(a,b)$
and column $k(a',b')$ is
\begin{itemize}
\item $1$, if $a'=a-1$, $b'=b-1$, or $a'=a+1$, $b'=b$, or $a'=a$, $b'=b+1$;
\item $-1$, if $a'=a+1$, $b'=b+1$, or $a'=a-1$, $b'=b$, or $a'=a$,
$b'=b-1$;
\item $0$, otherwise.
\end{itemize}

Let $B=(b_{ij})$\@. The mutation rule is given by:
$$
c_k'=\frac{\prod_{b_{ik}>0}c_i^{b_{ik}}
+\prod_{b_{ik}<0}c_i^{-b_{ik}}}{c_k}
$$
and if $\mu_k(B)=(b'_{ij})$ by
\begin{equation*}
\label{eq:mutation}
b'_{ij} =
\begin{cases}
-b_{ij} & \text{if $i=k$ or $j=k$;} \\[.05in]
b_{ij} + \displaystyle\frac{|b_{ik}| b_{kj} +
b_{ik} |b_{kj}|}{2} & \text{otherwise.}
\end{cases}
\end{equation*}

It is not at all obvious from the above definition that mutation is 
well-defined, namely why each~$c'_k$ is in $\palg$\@. However
Berenstein-Fomin-Zelevinsky \cite{BFZ} showed that this is true. In
particular, any sequence of mutations produces a cluster. Let us identify
two clusters if one can be produced from the other by simultaneous reindexing of 
the elements $c_k$ and the rows and columns of the matrix $B$\@. It is then
natural to ask whether the number of clusters is finite or infinite. It
turns out that it is finite for $n\leq 5$, but there are infinitely many
clusters for $n\geq 6$\@. 

The elements of the clusters are called {\em cluster variables}, while products
of the form $\prod_{c_k\in\cluster} c_k^{a_k}$ are called {\em cluster
monomials}.  These generate the algebra $\palg$.  
For $n\leq 5$, the cluster monomials span $\palg$ as a vector space, 
but for $n\geq 6$ they are expected to span a proper subspace.  

We have programmed this inductive definition into the computer, 
and used it to compute all cluster variables in types $\rm A_2, A_3, A_4$ 
(where there are $4,12,40$ of them, respectively)
and a list of~$719$ cluster variables in type $\rm A_5$~\cite{AK*}\@.

\end{document}